\documentclass[reqno,final,11pt]{amsart}
\usepackage{natbib}  
\usepackage{fancyhdr} 
\usepackage{color} 
\usepackage{hyperref} 
\usepackage{graphicx} 
\usepackage{geometry}
\usepackage[utf8]{inputenc}
\usepackage[T1]{fontenc}
\usepackage{verbatim}
\usepackage{lineno}

\usepackage{amsmath}
\usepackage{amsfonts}
\usepackage{hyperref}
\usepackage{xcolor}
\usepackage{graphicx}
\usepackage{enumerate}
\usepackage{enumitem}
\usepackage{xcolor}
\usepackage{tikz}
\usetikzlibrary{calc,shapes,backgrounds,arrows}
\usepackage{amssymb}

\definecolor{aleacolor}{rgb}{0.16,0.59,0.78}

\hypersetup{
breaklinks,
colorlinks=true,
linkcolor=aleacolor,
urlcolor=aleacolor,
citecolor=aleacolor}

\newcommand{\R}{\mathbb{R}}

\renewcommand{\cite}{\citet}

\theoremstyle{plain}
\newtheorem{theorem}{Theorem}[section]                                          
\newtheorem{proposition}[theorem]{Proposition}                          
\newtheorem{lemma}[theorem]{Lemma}

\theoremstyle{definition}

\theoremstyle{remark}
\newtheorem{remark}[theorem]{Remark}
\newtheorem{example}[theorem]{Example}

\makeatletter \@addtoreset{equation}{section} \makeatother
\renewcommand\theequation{\thesection.\arabic{equation}}
\renewcommand\thefigure{\thesection.\arabic{figure}}
\renewcommand\thetable{\thesection.\arabic{table}}

\def\R{\mathbb{R}}

\def\L{\mathcal{L}}

\newtheorem{condition}{Condition}[section]

\begin{document}

\title[Critical multitype Bellman--Harris processes with finite-mean lifetimes]
{Convergence of a Critical Multitype Bellman--Harris Process with Finite-Mean Lifetimes}

\author{E.T. Kolkovska}
\address{Centro de Investigaci\'{o}n en Matem\'{a}ticas, Guanajuato, Mexico.}
\email{todorova@cimat.mx} 

\author{J.A. L\'{o}pez-Mimbela}
\address{Centro de Investigaci\'{o}n en Matem\'{a}ticas, Guanajuato, Mexico.}
\email{jalfredo@cimat.mx} 

\author{J.H. Ram\'{\i}rez-Gonz\'{a}lez}
\address{Universidade de S\~ao Paulo, Brasil}
\email{hermenegildo.ramirez@usp.br}


\begin{abstract}
We study a critical multitype Bellman--Harris branching particle system in
\(\mathbb R^N\) with a finite type space \(\mathbf K=\{1,\dots,K\}\). Particles of type \(i\)
move according to a symmetric \(\alpha_i\)-stable process, have non-arithmetic lifetimes with
finite mean, and reproduce according to a critical offspring law whose mean matrix is
irreducible and stochastic. The branching mechanism is assumed to be in the domain of
attraction of a \((1+\beta)\)-stable law, with \(\beta\in(0,1]\). We prove that the particle
system converges, as \(t\to\infty\), to a limiting random measure which is nonzero. We also
show that the limiting population preserves the initial intensity measure. Thus, the system
persists with full intensity. These results complement the local extinction results of
\cite{Kevei} for critical multitype Bellman--Harris systems.
\end{abstract}

\maketitle

\begingroup
\renewcommand{\thefootnote}{}
\footnotetext[1]{Corresponding author: \texttt{hermenegildo.ramirez@usp.br}.}
\endgroup

\section{Introduction and previous work}

A central question for population models based on branching processes is the classification of
their long-time behaviour as random measures. In critical branching systems, where the mean
number of offspring per individual is one, different asymptotic regimes may occur and should be
kept distinct.

Let $E$ be a locally compact state space. A particle system on $E$ will be denoted by
$(X_t)_{t\ge0}$ and viewed as an $\mathcal M(E)$-valued stochastic process, where
$\mathcal M(E)$ is the space of nonnegative Radon measures on $E$, endowed with the vague
topology. For $\mu\in\mathcal M(E)$ and a measurable function $\varphi:E\to\mathbb R$, we write
\[
\langle \mu,\varphi\rangle:=\int_E \varphi\,d\mu .
\]
We also set
\[
C_c^+(E):=
\left\{
\varphi:E\to[0,\infty):\ \varphi \text{ is continuous and has compact support}
\right\}.
\]

We say that the particle system $(X_t)_{t\ge0}$ survives in the limiting sense if
\[
X_t \Rightarrow X_\infty,\qquad t\to\infty,
\]
in $\mathcal M(E)$, for some random measure $X_\infty\in\mathcal M(E)$ such that
\[
\mathbb P\{X_\infty\neq 0\}>0.
\]
If, in addition, the limiting random measure preserves the initial mean measure, in the sense that
\[
\mathbb E\langle X_\infty,\varphi\rangle
=
\mathbb E\langle X_0,\varphi\rangle,
\qquad \varphi\in C_c^+(E),
\]
then we say that the system persists with full intensity. Finally, we say that the system undergoes
local extinction if
\[
\langle X_t,\varphi\rangle
\xrightarrow[t\to\infty]{\mathbb P}0,
\qquad \varphi\in C_c^+(E).
\]

In \cite{GW}, a monotype system in $\mathbb R^N$ is studied in which particles perform
symmetric $\alpha$-stable motions, $0<\alpha\le2$, and branch at the end of exponential lifetimes.
The offspring law has generating function
\[
H(s)=s+\frac12(1-s)^{1+\beta},\qquad 0<\beta\le1.
\]
Starting from a homogeneous Poisson population, they prove local extinction when
$N\le \alpha/\beta$, and convergence to a nonzero limiting equilibrium of full intensity when
$N>\alpha/\beta$. Later, \cite{VW} considered general non-arithmetic lifetimes with finite mean
and obtained the same critical dimension. Moreover, for heavy-tailed lifetimes with tail exponent
$0<\gamma\le1$, they proved local extinction if $N<\alpha\gamma/\beta$, convergence to a
Poisson equilibrium of full intensity if $N>\alpha\gamma/\beta$, and a nonzero subsequential
limit in the critical case $N=\alpha\gamma/\beta$; see also \cite{FVW} for further results at the
critical dimension.

 We now turn to multitype systems. Consider a branching population in $\mathbb R^N$ with types
$i\in\mathbf K:=\{1,\dots,K\}$. A type-$i$ particle moves according to a symmetric
$\alpha_i$-stable motion, lives for a non-arithmetic time with distribution $F_i$, and then branches
according to the multitype generating function $f_i(\mathbf s)$, $\mathbf s\in[0,1]^K$. Let
$M=(m_{i,j})_{i,j=1}^K$ be the mean matrix, where
\[
m_{i,j}=\frac{\partial f_i}{\partial s_j}(\mathbf 1),
\qquad \mathbf 1=(1,\dots,1).
\]

In the Markovian case, where type-$i$ lifetimes are exponential with mean $1/\lambda_i$, set
$r_{i,j}:=\lambda_i m_{i,j}$ for $i\ne j$, and consider a Poisson initial field with intensity
\[
\bar\Lambda:=\sum_{i=1}^K \theta_i(\lambda\otimes\delta_i),
\]
where $\lambda$ is Lebesgue measure on $\mathbb R^N$. A standard requirement for
$\bar\Lambda$ to be a stationary intensity is that
$\theta=(\theta_i)_{i=1}^K$ be a stationary distribution for the type chain with jump rates
$(r_{i,j})_{i\ne j}$, equivalently,
\begin{equation}\label{estatio}
\sum_{i\ne j}\theta_i r_{i,j}
=
\theta_j\sum_{i\ne j} r_{j,i},
\qquad j=1,\dots,K.
\end{equation}

For $\phi\in C_c^+(\mathbb R^N\times\mathbf K)$ and
$(x,i)\in \mathbb R^N\times\mathbf K$, write
\[
W_t(x,i):=
E\!\left(
e^{-\langle X_t,\phi\rangle}
\,\middle|\, X_0=\delta_{(x,i)}
\right).
\]
If $X_0$ is Poisson with intensity $\bar \Lambda$, then \cite{GW,GRW} identify the Laplace
functional as
\begin{equation}\label{persistencea}
E\!\left(e^{-\langle X_t,\phi\rangle}\right)
=
\exp\{-\langle\bar\Lambda,1-W_t(\cdot)\rangle\}.
\end{equation}
Using this representation, together with asymptotic estimates and tightness, they obtain
convergence in distribution of $X_t$ to an equilibrium $X_\infty$. In \cite{Lopez1}, assuming that for each type $i$ the total offspring number belongs to the normal domain of attraction of a $(1+\beta_i)$-stable law ($0<\beta_i\le1$), they recover the same extinction/convergence dichotomy with critical threshold expressed through the parameters $\alpha:=\min_i\alpha_i$ and $\beta:=\min_i\beta_i$.

When lifetimes are not exponential, the process is no longer Markovian on $\R^N\times\mathbf K$ and one is led to multitype Bellman--Harris systems. In \cite{Kevei}, the mean matrix $M$ is assumed to be an ergodic stochastic matrix (hence critical, with Perron--Frobenius eigenvalue $1$), lifetimes are non-arithmetic with finite mean, and the offspring mechanism satisfies regular-variation assumptions. More precisely, if
\begin{equation}\label{kevei-cond_1}
x-\langle v,1-\mathbf f(1-ux)\rangle\sim x^{1+\beta}L(x)\ \ (x\to0)
\end{equation}
and
\begin{equation}\label{kevei-cond_2}
\lim_{n\to\infty}\frac{n(1-F_i(n))}{\langle v,1-\mathbf f^{(n)}(0)\rangle}=0\ \ (i=1,\dots,K),
\end{equation}
where \(u\) and \(v\) are respectively the normalized right and left eigenvectors of
\(M\) associated with the eigenvalue \(1\), \(L\) is slowly varying at \(0\), and
\(\beta\in(0,1]\), then local extinction holds whenever $N<\alpha/\beta$ with $\alpha:=\min_{1\le i\le K}\alpha_i$. Thus, under fairly general hypotheses, the extinction side is understood for multitype non-Markovian systems, while convergence to nonzero limiting equilibria is much less developed.

The purpose of the present paper is to address this convergence problem for a
critical multitype Bellman--Harris branching particle system with finite mean
lifetimes. We consider a finite type space \(\mathbf K=\{1,\dots,K\}\), type-dependent
symmetric stable motions, and a critical offspring mechanism whose mean matrix is
irreducible and stochastic. The branching mechanism is assumed to be in the domain
of attraction of a \((1+\beta)\)-stable law, in the sense of \eqref{kevei-cond_1}.
For a Poisson initial population with the stationary intensity weights determined by
the finite mean lifetimes, we prove that, under the Condition~\ref{Ac} and if
\[
q\wedge \frac{N}{\alpha_*}>\frac1\beta,
\qquad
\alpha_*:=\alpha_1\wedge\cdots\wedge\alpha_K,
\]
the system converges, as \(t\to\infty\), to a nonzero limiting random measure. We
also prove that the limiting population preserves the initial intensity measure, and
therefore the system persists with full intensity. Finally, we provide a sufficient
condition for the occupation-time estimate in terms of polynomial tail bounds for
the lifetime distributions.

\section{The model and main results}\label{results}

\subsection{The model}\label{subsec:model}

Let \(K\ge2\) and let
\[
\mathbf K:=\{1,\dots,K\},
\qquad
S:=\mathbb R^N\times\mathbf K.
\]
The population at time \(t\ge0\) is represented by the locally finite counting measure
\[
X_t=\sum_k\delta_{(X_t^k,K_t^k)}\in\mathcal M(S),
\]
where \(X_t^k\in\mathbb R^N\) is the position of the \(k\)-th particle alive at time \(t\), and
\(K_t^k\in\mathbf K\) is its type.

\medskip

\noindent\textbf{Motion.}
A particle of type \(i\in\mathbf K\) moves according to a symmetric
\(\alpha_i\)-stable process \((\xi_t^{(i)})_{t\ge0}\) in \(\mathbb R^N\), started from its birth
position. The motion indices satisfy
\[
\alpha_i\in(0,2],
\qquad i=1,\dots,K.
\]
The motions of distinct particles are independent, and are also independent of the lifetimes and
offspring variables.

\medskip

\noindent\textbf{Lifetimes.}
The lifetime of a type-\(i\) particle is denoted by \(\tau_i\) and has distribution function
\(F_i\) on \([0,\infty)\). We assume throughout that
\[
F_i(0)=0,
\qquad
F_i \text{ is non-arithmetic},
\qquad
m_i:=E(\tau_i)=\int_0^\infty t\,dF_i(t)<\infty,
\qquad i=1,\dots,K.
\]
Lifetimes are independent of each other and of all motions and offspring variables.

\medskip

\noindent\textbf{Branching mechanism.}
At the end of its lifetime, a type-\(i\) particle is replaced by an offspring vector
\[
\boldsymbol\zeta_i=(\zeta_{i,1},\dots,\zeta_{i,K})\in\mathbb N_0^K,
\]
where \(\zeta_{i,j}\) is the number of type-\(j\) children. Its generating function is
\[
f_i(\mathbf s)
:=
E\left(\prod_{j=1}^K s_j^{\zeta_{i,j}}\right),
\qquad
\mathbf s=(s_1,\dots,s_K)\in[0,1]^K.
\]
Let
\[
\mathbf f:=(f_1,\dots,f_K),
\]
and let \(M=(m_{i,j})_{i,j=1}^K\) be the mean offspring matrix, where
\[
m_{i,j}
:=
\frac{\partial f_i}{\partial s_j}(\mathbf 1),
\qquad
\mathbf 1=(1,\dots,1).
\]
We assume that \(M\) is irreducible and stochastic, and that
\[
m_{i,i}>0,
\qquad i=1,\dots,K.
\]
Thus the system is critical. Let \(\pi=(\pi_1,\dots,\pi_K)\) be the stationary distribution of \(M\):
\[
\pi M=\pi,
\qquad
\sum_{i=1}^K\pi_i=1,
\qquad
\pi_i>0,\quad i=1,\dots,K.
\]
Since \(M\) is stochastic, we set
\[
u:=\mathbf 1=(1,\dots,1),
\qquad
v:=\pi.
\]
Then
\[
Mu=u,
\qquad
vM=v,
\qquad
\langle v,u\rangle=1.
\]
We assume that, for some \(\beta\in(0,1]\) and some positive slowly varying function \(L\) at zero,
\begin{equation*}
x-\langle v,1-\mathbf f(1-u x)\rangle
\sim
x^{1+\beta}L(x),
\qquad x\downarrow0.
\end{equation*}

\medskip

\noindent\textbf{Initial configuration.}
We consider a Poisson initial population on \(S\). The weights of the initial intensity are chosen
according to the stationary occupation proportions of the type process in continuous time.

Set
\[
\bar m:=\sum_{i=1}^K\pi_i m_i,
\qquad
\theta_i:=\frac{\pi_i m_i}{\bar m},
\qquad i=1,\dots,K.
\]
Equivalently, these weights satisfy the balance relation obtained from the continuous-time rates
\[
r_{i,j}:=\frac{m_{i,j}}{m_i},
\qquad i\ne j.
\]
Indeed, for every \(j=1,\dots,K\),
\begin{equation*}
\sum_{i\ne j}\theta_i r_{i,j}
=
\theta_j\sum_{k\ne j}r_{j,k}.
\end{equation*}
To verify this, fix \(j\). Since \(\theta_i=\pi_i m_i/\bar m\),
\[
\sum_{i\ne j}\theta_i r_{i,j}
=
\sum_{i\ne j}
\frac{\pi_i m_i}{\bar m}
\frac{m_{i,j}}{m_i}
=
\frac1{\bar m}
\sum_{i\ne j}\pi_i m_{i,j}.
\]
Using \(\pi M=\pi\), we have
\[
\pi_j
=
\sum_{i=1}^K\pi_i m_{i,j}
=
\pi_jm_{j,j}
+
\sum_{i\ne j}\pi_i m_{i,j}.
\]
Hence
\[
\sum_{i\ne j}\pi_i m_{i,j}
=
\pi_j(1-m_{j,j})
=
\pi_j\sum_{k\ne j}m_{j,k}.
\]
Therefore,
\[
\sum_{i\ne j}\theta_i r_{i,j}
=
\frac{\pi_j}{\bar m}\sum_{k\ne j}m_{j,k}.
\]
On the other hand,
\[
\theta_j\sum_{k\ne j}r_{j,k}
=
\frac{\pi_jm_j}{\bar m}
\sum_{k\ne j}\frac{m_{j,k}}{m_j}
=
\frac{\pi_j}{\bar m}\sum_{k\ne j}m_{j,k}.
\]
This proves \eqref{estatio}.

Thus, throughout the finite-mean setting, we take \(X_0\) to be a Poisson random
measure with intensity
\[
\bar\Lambda
:=
\sum_{i=1}^K\theta_i\,\lambda\otimes\delta_i,
\]
where \(\lambda\) denotes Lebesgue measure on \(\mathbb R^N\). The initial population is
independent of all motions, lifetimes, and offspring variables.

\medskip

We shall use the following quantitative condition on the associated type process. Let
\(Z=(Z_t)_{t\ge0}\) be the Markov renewal process on \(\mathbf K\) such that, when \(Z\) is in
state \(i\), the holding time has distribution \(F_i\), and after this holding time the process jumps
from \(i\) to \(j\) with probability \(m_{i,j}\). Let
\[
L_i(t):=\int_0^t\mathbf 1_{\{Z_s=i\}}\,ds
\]
be the time spent by \(Z\) in state \(i\) up to time \(t\).

\begin{condition}\label{Ac}
For each \(i\in\mathbf K\), there exist constants \(c_i\in(0,1)\), \(C>0\), and \(q>0\) such that,
for all sufficiently large \(t\) and every initial type \(j\in\mathbf K\),
\begin{equation}\label{eq:occupation-condition}
P_j\left(\frac{L_i(t)}{t}\le c_i\right)
\le
C t^{-q},
\end{equation}
where \(P_j\) denotes the law of \(Z\) started from \(j\).
\end{condition}

For \(\phi\in C_c^+(S)\), define
\[
W_t(x,i):=E_{(x,i)}\left(e^{-\langle X_t,\phi\rangle}\right),
\qquad
V_t(x,i):=1-W_t(x,i),
\qquad (x,i)\in S.
\]
We also write
\[
\|V_t(\cdot,i)\|_1:=\int_{\mathbb R^N}V_t(x,i)\,dx,
\qquad
G_i(t):=\|V_t(\cdot,i)\|_1,
\qquad
g_\phi(t):=\sum_{i=1}^K\theta_iG_i(t).
\]

\medskip

\subsection{Main results}\label{subsec:main-results}

Throughout this subsection, the lifetime distributions are assumed to be in the finite-mean
regime specified in Section~\ref{subsec:model}: for each $i=1,...,K$, $F_i$ is non-arithmetic and satisfies
$m_i:=E(\tau_i)<\infty$.

\begin{proposition}[Finite-mean first-moment bound and Poisson Laplace functional]
\label{Lema2.6.1}
Consider the Bellman--Harris branching particle system described in
Section~\ref{subsec:model}, with Poisson initial intensity \(\bar\Lambda\). For
\(\phi\in C_c^+(S)\), let \(V_t\), \(G_i(t)\), and \(g_\phi(t)\) be defined as in
Section~\ref{subsec:model}. Then the following properties hold.

\begin{itemize}
\item[\textup{(a)}]  For every \(\phi\in C_c^+(S)\),
 \[
 \sup_{t\ge0}
 \sum_{i=1}^K
 \int_{\mathbb R^N}
 E_{(x,i)}\big(\langle X_t,\phi\rangle\big)\,dx
 <\infty.
 \]
 Consequently, for each \(i=1,\ldots,K\), the function
 \(G_i(t)=\|V_t(\cdot,i)\|_1\) is bounded on \([0,\infty)\). Moreover, for every
 \(\phi\in C_c^+(S)\),
 \begin{equation}\label{elimit}
 E\big(\langle X_t,\phi\rangle\big)
 \longrightarrow
 \sum_{i=1}^K\theta_i\|\phi_i\|_1,
 \qquad t\to\infty .
 \end{equation}
In particular, \((X_t)_{t\ge0}\) is tight in the space of Borel--Radon measures on
\(S\), endowed with the vague topology.

\item[\textup{(b)}]
For every \(\phi\in C_c^+(S)\),
\begin{equation}\label{laplace}
\begin{split}
E\left(e^{-\langle X_t,\phi\rangle}\right)
&=
\exp\left\{
-\sum_{i=1}^K\theta_i
\int_{\mathbb R^N}
E_{(x,i)}
\left(
1-e^{-\langle X_t,\phi\rangle}
\right)\,dx
\right\}
\\
&=
\exp\left\{
-\sum_{i=1}^K\theta_i\|V_t(\cdot,i)\|_1
\right\}
=
e^{-g_\phi(t)}.
\end{split}
\end{equation}
\end{itemize}
\end{proposition}

The following result gives a simple sufficient condition for Condition~\ref{Ac} in terms of polynomial
tail bounds for the lifetime distributions.

\begin{theorem}\label{Teo2.3}
Assume the model described in Section~\ref{subsec:model}. Assume moreover that the
    lifetime distributions satisfy
\begin{equation}\label{tail2}
F_i(0)=0,\qquad
F_i(t)<1 \text{ for all } t\ge0,\qquad
1-F_i(t)\le A t^{-\eta^{(i)}} \text{ for all } t>0,
\end{equation}
for some \(A>0\) and some \(\eta^{(i)}>1\), \(i=1,\dots,K\). Set
\[
\eta:=\min_{1\le i\le K}\eta^{(i)}.
\]
Then Condition~\ref{Ac} holds with \(q=\eta-1\).
\end{theorem}

Theorem~\ref{Teo2.3} shows that Condition~\ref{Ac} holds for a broad class of finite-mean lifetime
distributions. 

\begin{theorem}\label{Teorema2.6.8}
Consider the Bellman--Harris branching particle system described in
Section~\ref{subsec:model}. Assume that Condition~\ref{Ac} holds for some
\(q>0\).

Assume that the branching mechanism satisfies
\begin{equation}\label{eq:K-kevei-cond}
x-\langle v,1-\mathbf f(1-u x)\rangle
\sim
x^{1+\beta}L(x),
\qquad x\downarrow0,
\end{equation}
for some \(\beta\in(0,1]\) and some positive slowly varying function \(L\) at \(0\). Let
\[
\alpha_*:=\alpha_1\wedge\cdots\wedge\alpha_K,
\qquad
\rho:=q\wedge\frac{N}{\alpha_*}.
\]
Assume that
\[
\rho>\frac1\beta.
\]

For \(\phi\in C_c^+(S)\), let \(V_t\), \(G_i(t)\), and \(g_\phi(t)\) be defined as in
Section~\ref{subsec:model}. For \(i=1,\dots,K\), define
\[
R_i(\mathbf y):=
\sum_{j=1}^K m_{i,j}y_j
-
\bigl(1-f_i(1-\mathbf y)\bigr),
\qquad \mathbf y\in[0,1]^K,
\]
and
\[
B_i(t):=
\int_{\mathbb R^N}
R_i\bigl(V_t(x,1),\dots,V_t(x,K)\bigr)\,dx.
\]
Then, for every nonzero \(\phi\in C_c^+(S)\), the limit
\begin{equation}\label{gfinal}
G(\phi):=\lim_{t\to\infty}g_\phi(t)
=
\sum_{k=1}^K\theta_kG_k(0)
-
\frac1{\bar m}
\sum_{k=1}^K\pi_k\int_0^\infty B_k(s)\,ds
\end{equation}
exists and satisfies \(G(\phi)>0\).

Moreover, if
\[
\phi_\epsilon:=-\log\bigl(1-\epsilon(1-e^{-\phi})\bigr),
\qquad 0<\epsilon<1,
\]
then
\begin{equation}\label{eq:K-first-order-limit}
\lim_{\epsilon\downarrow0}
\frac{G(\phi_\epsilon)}{\epsilon}
=
\sum_{i=1}^K\theta_i\|1-e^{-\phi_i}\|_1.
\end{equation}
\end{theorem}

Recall that \(S=\mathbb{R}^N\times\mathbb \mathbf K\), where
\(\mathbf K=\{1,\ldots,K\}\) is finite. Endowing \(\mathbf K\) with the
discrete topology, it follows that \(S\) is a locally compact second
countable Hausdorff space. Hence the convergence criterion for random
measures in the vague topology, Theorem~4.2 of \cite{Olav}, applies on
\(S\).

By Proposition~\ref{Lema2.6.1} and Theorem~\ref{Teorema2.6.8}, for every \(\phi\in C_c^+(S)\),
\[
  \mathbb E\exp\{-\langle X_t,\phi\rangle\}
  \longrightarrow
  \mathbb E\exp\{-\langle X_\infty,\phi\rangle\},
\]
where \(X_\infty\) is a random Radon measure on \(S\) whose Laplace
functional is given by
\[
  \mathbb E\big(e^{-\langle X_\infty,\phi\rangle}\big)
  = e^{-G(\phi)}.
\]
Thus the Laplace functionals converge on \(C_c^+(S)\), the class of
nonnegative continuous functions with compact support. Therefore, by
Theorem~4.2 of \cite{Olav},
\[
  X_t \Rightarrow X_\infty, \qquad t\to\infty,
\]
in the space of Radon measures on \(S\) endowed with the vague topology.

\begin{remark}
Let \(X_\infty\) be the limiting random measure in
Theorem~\ref{Teorema2.6.8}, characterized by
\[
E\bigl(e^{-\langle X_\infty,\phi\rangle}\bigr)
=
e^{-G(\phi)},
\qquad
\phi\in C_c^+(S).
\]

The limiting population preserves the initial intensity measure. Indeed, let \(\psi\in C_c^+(S)\). Choose \(a>0\) such that
\[
0\le a\psi<1,
\]
and define
\[
\phi_a:=-\log(1-a\psi).
\]
Then
\[
1-e^{-\phi_a}=a\psi.
\]
For \(0<\varepsilon<1\), set
\[
\phi_{a,\varepsilon}
:=
-\log\bigl(1-\varepsilon(1-e^{-\phi_a})\bigr)
=
-\log(1-\varepsilon a\psi).
\]
By \eqref{eq:K-first-order-limit},
\[
\lim_{\varepsilon\downarrow0}
\frac{G(\phi_{a,\varepsilon})}{\varepsilon}
=
\sum_{i=1}^K\theta_i\|a\psi_i\|_1
=
a\int_S\psi\,d\bar\Lambda .
\]
Since \(G(\phi_{a,\varepsilon})\to0\), it follows that
\[
\lim_{\varepsilon\downarrow0}
\frac{1-e^{-G(\phi_{a,\varepsilon})}}{\varepsilon}
=
a\int_S\psi\,d\bar\Lambda .
\]
On the other hand, by the Laplace functional of \(X_\infty\),
\[
1-e^{-G(\phi_{a,\varepsilon})}
=
E\left(
1-e^{-\langle X_\infty,\phi_{a,\varepsilon}\rangle}
\right).
\]
Therefore,
\begin{equation}\label{eq:intensity-limit-remark}
\lim_{\varepsilon\downarrow0}
E\left[
\frac{
1-e^{-\langle X_\infty,\phi_{a,\varepsilon}\rangle}
}{\varepsilon}
\right]
=
a\int_S\psi\,d\bar\Lambda .
\end{equation}

For each fixed Radon measure \(\mu\in\mathcal M(S)\),
\[
\frac{\phi_{a,\varepsilon}}{\varepsilon}
=
\frac{-\log(1-\varepsilon a\psi)}{\varepsilon}
\longrightarrow
a\psi
\]
pointwise, and since \(\psi\) has compact support,
\[
\frac{
1-e^{-\langle \mu,\phi_{a,\varepsilon}\rangle}
}{\varepsilon}
\longrightarrow
a\langle \mu,\psi\rangle .
\]
Applying this with \(\mu=X_\infty\), Fatou's lemma and
\eqref{eq:intensity-limit-remark} give
\[
aE\langle X_\infty,\psi\rangle
\le
a\int_S\psi\,d\bar\Lambda
<\infty.
\]
Hence \(E\langle X_\infty,\psi\rangle<\infty\).

Moreover, since \(a\|\psi\|_\infty<1\), there exists \(C_a>0\) such that, for all
\(0<\varepsilon<1\),
\[
\frac{\phi_{a,\varepsilon}}{\varepsilon}
\le
C_a\,a\psi .
\]
Using \(1-e^{-x}\le x\), we obtain
\[
0\le
\frac{
1-e^{-\langle X_\infty,\phi_{a,\varepsilon}\rangle}
}{\varepsilon}
\le
\left\langle X_\infty,\frac{\phi_{a,\varepsilon}}{\varepsilon}\right\rangle
\le
C_a a\langle X_\infty,\psi\rangle .
\]
The right-hand side is integrable. Thus, by dominated convergence,
\[
\lim_{\varepsilon\downarrow0}
E\left[
\frac{
1-e^{-\langle X_\infty,\phi_{a,\varepsilon}\rangle}
}{\varepsilon}
\right]
=
aE\langle X_\infty,\psi\rangle .
\]
Combining this identity with \eqref{eq:intensity-limit-remark}, we obtain
\[
aE\langle X_\infty,\psi\rangle
=
a\int_S\psi\,d\bar\Lambda .
\]
Since \(a>0\),
\[
E\langle X_\infty,\psi\rangle
=
\int_S\psi\,d\bar\Lambda .
\]
Thus the limiting population preserves the initial intensity measure:
\[
E(X_\infty)=\bar\Lambda .
\]
\end{remark}

\begin{example}\label{Example1}
Consider the case \(K=2\). Let the type-\(1\) lifetime distribution be
\(F_1:=F_{p,\lambda_1,\lambda_2}\), where \(F_{p,\lambda_1,\lambda_2}\) has density
\[
f_{p,\lambda_1,\lambda_2}(x)
:=
\Big(p\lambda_1e^{-\lambda_1 x}+(1-p)\lambda_2e^{-\lambda_2 x}\Big)
\mathbf 1_{(0,\infty)}(x),
\qquad
p\in[0,1],\quad \lambda_1,\lambda_2>0.
\]
Then
\begin{equation}\label{fcond1}
1-F_1(x)
=
pe^{-\lambda_1 x}+(1-p)e^{-\lambda_2 x}
\le
e^{-(\lambda_1\wedge\lambda_2)x},
\qquad x>0.
\end{equation}
Let the type-\(2\) lifetime distribution be Weibull with parameters \(r>0\) and \(\lambda>0\), with distribution function
\[
F_2(x)
:=
\left(1-e^{-(x/\lambda)^r}\right)\mathbf 1_{[0,\infty)}(x).
\]
Then
\begin{equation}\label{fcond2}
1-F_2(x)
=
e^{-(x/\lambda)^r},
\qquad x>0.
\end{equation}
It follows from \eqref{fcond1}--\eqref{fcond2} that, for every \(\eta>1\), there exists
\(A_\eta>0\) such that
\[
1-F_i(x)\le A_\eta x^{-\eta},
\qquad x>0,\quad i=1,2.
\]
Indeed, the exponential tail in \eqref{fcond1} and the Weibull tail in \eqref{fcond2}
decay faster than any negative power. Thus \eqref{tail2} holds for every \(\eta>1\), and
Theorem~\ref{Teo2.3} implies that Condition~\ref{Ac} holds with arbitrarily large \(q\).
Hence, if the branching mechanism satisfies \eqref{eq:K-kevei-cond} with exponent
\(\beta\in(0,1]\), then Theorem~\ref{Teorema2.6.8} gives persistence with full intensity whenever
\[
N>\frac{\alpha_1\wedge\alpha_2}{\beta}.
\]

For the extinction side, suppose in addition that the offspring mechanism belongs to the particular
class considered in Remark~1 of \cite{Kevei}, namely
\[
f_i(s,s)=s+c_i(1-s)^{1+\beta},
\qquad 0\le s\le1,\quad i=1,2,
\]
where \(c_i\in(0,(1+\beta)^{-1}]\).
Moreover, \eqref{fcond1}--\eqref{fcond2} imply that, for every \(\varepsilon>0\),
\[
\lim_{n\to\infty}
n^{1+\frac1\beta+\varepsilon}\bigl[1-F_i(n)\bigr]
=
0,
\qquad i=1,2.
\]
Therefore, by Remark~1 of \cite{Kevei}, condition \eqref{kevei-cond_2} is satisfied, and
Theorem~1 of \cite{Kevei} gives local extinction whenever
\[
N<\frac{\alpha_1\wedge\alpha_2}{\beta}.
\]
Consequently, for these two lifetime distributions and for the above particular offspring mechanism,
the extinction--persistence dichotomy is determined outside the critical dimension
\(N=\frac{\alpha_1\wedge\alpha_2}{\beta}\), which is not treated here.
\end{example}

\section{Proofs of main results}\label{PROOFS}

\subsubsection*{Proof of Proposition \ref{Lema2.6.1}:}
We first prove \textup{(a)}. Fix \(\phi\in C_c^+(S)\). Suppose that the
system starts at time \(0\) with one particle of type \(i\) at position
\(x\in\mathbb R^N\). A first-step decomposition gives
\begin{equation}\label{eq:first-moment-first-step}
\begin{split}
E_{(x,i)}\big(\langle X_t,\phi\rangle\big)
&=
(1-F_i(t))P_{t,\alpha_i}(\phi_i)(x)
\\
&\quad+
\int_{[0,t]}
P_{s,\alpha_i}
\left(
\sum_{j=1}^K
m_{i,j}
E_{(\cdot,j)}
\big(\langle X_{t-s},\phi\rangle\big)
\right)(x)
\,F_i(ds),
\end{split}
\end{equation}
where \(\phi_i(x):=\phi(x,i)\). This identity only uses the mean offspring
matrix \(M\).

Set
\[
H_i(t):=
\int_{\mathbb R^N}
E_{(x,i)}
\big(\langle X_t,\phi\rangle\big)\,dx,
\qquad i=1,\dots,K.
\]
Integrating \eqref{eq:first-moment-first-step} over \(\mathbb R^N\), and using
the invariance of Lebesgue measure under the \(\alpha_i\)-stable semigroup,
gives
\begin{equation}\label{eq:first-moment-integrated-system}
H_i(t)
=
(1-F_i(t))\|\phi_i\|_1
+
\sum_{j=1}^K
m_{i,j}
\int_{[0,t]}
H_j(t-s)\,F_i(ds),
\qquad i=1,\dots,K.
\end{equation}
Equivalently,
\[
H_i(t)
=
z_i^\phi(t)
+
\sum_{j=1}^K
\int_{[0,t]}H_j(t-s)\,F_{i,j}(ds),
\]
where
\[
z_i^\phi(t):=(1-F_i(t))\|\phi_i\|_1,
\qquad
F_{i,j}(t):=m_{i,j}F_i(t).
\]

We first justify local boundedness on compact time intervals. This is done
without assuming continuity of the lifetime distribution functions. Since
each distribution function is right-continuous and \(F_i(0)=0\), there exists
\(h>0\) such that
\[
\kappa_h:=\max_{1\le i\le K}F_i(h)<1.
\]
Let \(\mathcal B_h\) be the space of bounded Borel maps
\[
\mathbf u:[0,h]\to \mathbb R_+^K,
\]
endowed with the norm
\[
\|\mathbf u\|_{\infty,h}
:=
\sup_{0\le t\le h}\max_{1\le i\le K}|u_i(t)|.
\]
This is a complete metric space. Define, for
\(\mathbf u=(u_1,\ldots,u_K)\in\mathcal B_h\),
\[
(\mathcal T\mathbf u)_i(t)
:=
z_i^\phi(t)
+
\sum_{j=1}^K m_{i,j}
\int_{[0,t]}u_j(t-s)\,F_i(ds),
\qquad 0\le t\le h.
\]
The map \(\mathcal T\mathbf u\) is Borel measurable. Indeed, after extending
\(u_j\) by zero outside \([0,h]\), the function
\[
(t,s)\mapsto u_j(t-s)\mathbf 1_{\{0\le s\le t\}}
\]
is Borel measurable on \([0,h]\times[0,h]\), and \(F_i\) is a finite measure
on compact intervals. Since \(z_i^\phi\) is bounded on \([0,h]\), the operator
\(\mathcal T\) maps \(\mathcal B_h\) into itself. Moreover, for
\(\mathbf u,\mathbf v\in\mathcal B_h\),
\[
\begin{aligned}
\| \mathcal T\mathbf u-\mathcal T\mathbf v\|_{\infty,h}
&\le
\max_{1\le i\le K}
\sup_{0\le t\le h}
\sum_{j=1}^K m_{i,j}
\int_{[0,t]} |u_j(t-s)-v_j(t-s)|\,F_i(ds) \\
&\le
\max_{1\le i\le K}F_i(h)
\|\mathbf u-\mathbf v\|_{\infty,h} \\
&=
\kappa_h
\|\mathbf u-\mathbf v\|_{\infty,h}.
\end{aligned}
\]
Thus \(\mathcal T\) is a contraction on \(\mathcal B_h\).

To connect this contraction argument with the first moment \(H_i\), let
\[
\mathbf H^{[0]}\equiv 0,
\]
and define recursively
\[
\mathbf H^{[n+1]}:=\mathcal T\mathbf H^{[n]},
\qquad n\ge0.
\]
Equivalently,
\[
H_i^{[n+1]}(t)
=
z_i^\phi(t)
+
\sum_{j=1}^K m_{i,j}
\int_{[0,t]}H_j^{[n]}(t-s)\,F_i(ds),
\qquad 0\le t\le h.
\]
The sequence \((\mathbf H^{[n]})_{n\ge0}\) is the increasing sequence of
integrated first moments obtained from finite-depth truncations of the
branching construction. Hence, by monotone convergence,
\[
H_i^{[n]}(t)\uparrow H_i(t),
\qquad
0\le t\le h,\quad i=1,\ldots,K.
\]
On the other hand, since \(\mathcal T\) is a contraction on \(\mathcal B_h\),
 the iterates \(\mathbf H^{[n]}\) converge uniformly on \([0,h]\) to the
 unique fixed point \(\mathbf H^{[*]}\in\mathcal B_h\). Therefore
\[
H_i(t)=H_i^{[*]}(t),
\qquad 0\le t\le h,\quad i=1,\ldots,K,
\]
and consequently
\[
\sup_{0\le t\le h}\sum_{i=1}^K H_i(t)<\infty.
\]

Repeating the argument on consecutive intervals gives local boundedness on
every compact interval. Indeed, suppose that \(H_i\) is bounded on
\([0,a]\), for every \(i=1,\ldots,K\). For \(r\in[0,h]\), put
\[
Y_i(r):=H_i(a+r).
\]
From \eqref{eq:first-moment-integrated-system}, for \(0\le r\le h\),
\[
Y_i(r)
=
q_i^{(a)}(r)
+
\sum_{j=1}^K m_{i,j}
\int_{[0,r]}Y_j(r-s)\,F_i(ds),
\]
where
\[
q_i^{(a)}(r)
:=
z_i^\phi(a+r)
+
\sum_{j=1}^K m_{i,j}
\int_{(r,a+r]}H_j(a+r-s)\,F_i(ds).
\]
The function \(q_i^{(a)}\) is bounded on \([0,h]\), because
\(z_i^\phi\) is bounded on compact intervals and, for \(s\in(r,a+r]\), one
has \(a+r-s\in[0,a)\), where the functions \(H_j\) are bounded by the
induction hypothesis. Therefore, on \([0,h]\), the unknown vector
\(\mathbf Y=(Y_1,\ldots,Y_K)\) satisfies a renewal equation with bounded
Borel inhomogeneous term \(q^{(a)}\) and with the same contraction constant
\[
\max_{1\le i\le K}F_i(h)=\kappa_h<1.
\]

Starting from \(\mathbf Y^{[0]}\equiv0\) and iterating this equation, the
iterates converge uniformly on \([0,h]\) to a bounded fixed point, because the
contraction constant is again \(\kappa_h<1\). These iterates are precisely the first moments obtained from the corresponding
finite-depth truncations. Hence, by monotone convergence,
\[
Y_i^{[n]}(r)\uparrow Y_i(r),
\qquad 0\le r\le h,\quad i=1,\ldots,K.
\]
Thus \(\mathbf Y\) coincides with the bounded fixed point. Therefore
\[
\sup_{0\le r\le h}\sum_{i=1}^K Y_i(r)<\infty.
\]
Equivalently,
\[
\sup_{a\le t\le a+h}\sum_{i=1}^K H_i(t)<\infty.
\]
After finitely many steps, the functions \(H_i\) are bounded on any prescribed
compact interval.

We now apply the renewal theorem. The functions \(z_i^\phi\) are Borel
measurable, locally bounded, and directly Riemann integrable, since
\(1-F_i\) is nonincreasing and
\[
\int_0^\infty (1-F_i(t))\,dt=m_i<\infty.
\]
Moreover,
\[
F(\infty)=M,
\qquad
p(F(\infty))=p(M)=1,
\qquad
F(0)=0,
\qquad
p(F(0))=0<1.
\]
Since \(m_{i,i}>0\) and \(F_i\) is non-arithmetic, the diagonal measure
\(F_{i,i}=m_{i,i}F_i\) is non-arithmetic. Hence the matrix of measures
\((F_{i,j})_{i,j=1}^K\) is non-lattice. Therefore,
Theorem~2.2\textup{(iii)} of \cite{Athreya2} applies to
\eqref{eq:first-moment-integrated-system}.

Since \(M\) is stochastic, its right Perron--Frobenius eigenvector associated
with the eigenvalue \(1\) is
\[
u=\mathbf 1=(1,\dots,1),
\]
and its normalized left Perron--Frobenius eigenvector is
\[
v=\pi,
\qquad
\langle v,u\rangle=\sum_{i=1}^K\pi_i=1.
\]
Furthermore,
\[
b_{k,r}
:=
\int_0^\infty t\,dF_{k,r}(t)
=
m_{k,r}m_k.
\]
Thus the denominator appearing in the renewal theorem is
\[
\sum_{k=1}^K\sum_{r=1}^K v_k u_r b_{k,r}
=
\sum_{k=1}^K\sum_{r=1}^K \pi_k m_{k,r}m_k
=
\sum_{k=1}^K\pi_km_k
=
\bar m,
\]
because \(\sum_{r=1}^K m_{k,r}=1\). Also,
\[
\int_0^\infty z_k^\phi(s)\,ds
=
\|\phi_k\|_1
\int_0^\infty (1-F_k(s))\,ds
=
m_k\|\phi_k\|_1.
\]
Consequently, Theorem~2.2\textup{(iii)} of \cite{Athreya2} yields, for every
\(i=1,\dots,K\),
\begin{equation}\label{principal}
H_i(t)
\longrightarrow
\frac{\sum_{k=1}^K\pi_km_k\|\phi_k\|_1}{\bar m}
=
\sum_{k=1}^K\theta_k\|\phi_k\|_1,
\qquad t\to\infty.
\end{equation}
Combining this convergence with the local boundedness proved above, we obtain
\begin{equation}\label{eq:first-moment-uniform-bound}
\sup_{t\ge0}\sum_{i=1}^K H_i(t)<\infty.
\end{equation}

Since \(X_0\) is Poisson with intensity \(\bar\Lambda\),
\[
E\big(\langle X_t,\phi\rangle\big)
=
\sum_{i=1}^K\theta_i H_i(t).
\]
Using \eqref{principal} and \(\sum_{i=1}^K\theta_i=1\), we obtain
\[
E\big(\langle X_t,\phi\rangle\big)
\longrightarrow
\sum_{k=1}^K\theta_k\|\phi_k\|_1,
\qquad t\to\infty,
\]
which proves \eqref{elimit}.

Moreover,
\begin{equation}\label{aux1}
0\le
V_t(x,i)
=
E_{(x,i)}
\left(
1-e^{-\langle X_t,\phi\rangle}
\right)
\le
E_{(x,i)}
\big(\langle X_t,\phi\rangle\big).
\end{equation}
Integrating \eqref{aux1} over \(\mathbb R^N\), and using
\eqref{eq:first-moment-uniform-bound}, gives
\begin{equation}\label{eq:V-uniform-L1-bound}
\sup_{t\ge0}\sum_{i=1}^K\|V_t(\cdot,i)\|_1<\infty .
\end{equation}
Therefore each \(G_i(t)=\|V_t(\cdot,i)\|_1\) is bounded on \([0,\infty)\).

It remains to prove tightness. Let \(\mathcal K\subset S\) be compact, and
choose \(\psi\in C_c^+(S)\) such that \(\psi\ge 1_{\mathcal K}\). Then, by
Markov's inequality,
\[
P\big(X_t(\mathcal K)>a\big)
\le
P\big(\langle X_t,\psi\rangle>a\big)
\le
\frac{E(\langle X_t,\psi\rangle)}{a}.
\]
Since \(E(\langle X_t,\psi\rangle)\) is bounded uniformly in \(t\), it follows
that
\[
\lim_{a\to\infty}\limsup_{t\to\infty}
P\big(X_t(\mathcal K)>a\big)=0.
\]
Thus, by Lemma~4.5 of \cite{Olav}, the family \((X_t)_{t\ge0}\) is tight in
the vague topology. This proves \textup{(a)}.

We now prove \textup{(b)}. Similarly to \cite[page~622]{VW}, which in turn
relies on Remark~1.4.7 in \cite{LMW1988}, since the initial population
\(X_0\) is a Poisson population with intensity \(\bar\Lambda\) and is
independent of the branching mechanism, the Poisson exponential formula gives
\[
\begin{split}
E\left(e^{-\langle X_t,\phi\rangle}\right)
&=
\exp\left\{
-\sum_{i=1}^K\theta_i
\int_{\mathbb R^N}
E_{(x,i)}
\left(
1-e^{-\langle X_t,\phi\rangle}
\right)\,dx
\right\}
\\
&=
\exp\left\{
-\sum_{i=1}^K\theta_i\|V_t(\cdot,i)\|_1
\right\}.
\end{split}
\]
This proves \eqref{laplace}.
\(\hfill\square\)
\subsubsection*{Proof of Theorem \ref{Teo2.3}:}
Let \((Y_n)_{n\ge0}\) be the embedded Markov chain of the Markov renewal process
\((Z_t)_{t\ge0}\), with transition matrix \(M\), and with \(Y_0\) equal to the
initial type. Let \(\pi=(\pi_1,\dots,\pi_K)\) denote its stationary
distribution. Since \(M\) is irreducible and stochastic, \(\pi_j>0\) for every
\(j=1,\dots,K\). For \(j\in\{1,\dots,K\}\), define
\[
t_j(n):=\sum_{m=0}^{n-1}\mathbf 1_{\{Y_m=j\}},
\qquad n\ge1.
\]
By the large deviation theorem for finite Markov chains, see Lemma~2.13 of
\cite{Golds}, for every \(\epsilon>0\) there exist constants \(C,c>0\) such that
\begin{equation}\label{ldmc}
P_i\left(
\left|
\frac{t_j(n)}{n}-\pi_j
\right|>\epsilon
\right)
\le
Ce^{-cn},
\qquad n\ge1,\quad i,j\in\{1,\dots,K\}.
\end{equation}

Conditionally on the embedded chain \((Y_n)_{n\ge0}\), let
\((S_n)_{n\ge0}\) be independent holding times such that \(S_n\) has
distribution \(F_{Y_n}\). Define
\[
\sigma_0:=0,\qquad
\sigma_n:=\sum_{r=0}^{n-1}S_r,\quad n\ge1,
\]
and let
\[
n_t:=\max\{n\ge0:\sigma_n\le t\}
\]
be the number of renewals up to time \(t\). Then the occupation time of state
\(j\) up to time \(t\) is
\[
L_j(t)
=
\sum_{r=0}^{n_t-1}S_r\mathbf 1_{\{Y_r=j\}}
+
(t-\sigma_{n_t})\mathbf 1_{\{Y_{n_t}=j\}}.
\]
Equivalently, if \(T_1^j,T_2^j,\ldots\) denotes the sequence of holding times
attached to visits of the embedded chain to state \(j\), then
\[
L_j(t)
=
\sum_{k=1}^{t_j(n_t)}T_k^j
+
r_j(t),
\]
where
\[
r_j(t):=(t-\sigma_{n_t})\mathbf 1_{\{Y_{n_t}=j\}}\ge0.
\]
By construction, \(T_1^j,T_2^j,\ldots\) are i.i.d. with distribution \(F_j\).
Moreover, they are independent of the visit-count process
\((t_j(n))_{n\ge1}\), because the holding times are attached independently to
the visits of the embedded chain.

Fix \(j\in\{1,\dots,K\}\). Choose \(\epsilon>0\) such that
\[
\pi_j-\epsilon>0.
\]
Let
\[
\eta:=\min_{1\le r\le K}\eta^{(r)}.
\]
Let \(T^1,\ldots,T^K\) be independent random variables such that \(T^r\) has
distribution \(F_r\), and set
\[
T^*:=\max\{T^1,\dots,T^K\},
\qquad
\bar\mu:=E(T^*)<\infty.
\]
The finiteness of \(\bar\mu\) follows from Lemma~\ref{lemaaux1}, since
\(\eta>1\).

Choose \(c_{(j)}\in(0,1)\) sufficiently small so that
\[
d_j:=
\frac{2c_{(j)}}{m_j(\pi_j-\epsilon)}
<
\frac{1}{\bar\mu}.
\]
Equivalently, it is enough to take
\[
0<c_{(j)}
<
1\wedge
\frac{m_j(\pi_j-\epsilon)}{2\bar\mu}.
\]

We decompose
\[
\begin{split}
P_i\big(L_j(t)\le tc_{(j)}\big)
&\le
P_i\big(L_j(t)\le tc_{(j)},\, n_t>d_jt\big)
+
P_i\big(n_t\le d_jt\big).
\end{split}
\]

We first estimate the first term. By \eqref{ldmc},
\[
\begin{split}
P_i\left(
\frac{t_j(n_t)}{n_t}<\pi_j-\epsilon,\,
n_t>d_jt
\right)
&\le
\sum_{n\ge\lfloor d_jt\rfloor}
P_i\left(
\frac{t_j(n)}{n}<\pi_j-\epsilon
\right)
\\
&\le
\sum_{n\ge\lfloor d_jt\rfloor}Ce^{-cn}
\le
C_1e^{-c_1t}
\end{split}
\]
for all sufficiently large \(t\).

Next, we estimate the probability that the average of the type-\(j\) lifetimes is too
small. Since \(m_j=E(T_1^j)>0\), choose \(R>0\) such that
\[
E(T_1^j\wedge R)>\frac{3m_j}{4}.
\]
Set \(Y_k:=T_k^j\wedge R\). Then \(0\le Y_k\le R\), the variables
\(Y_k\) are independent, and
\[
E(Y_k)>\frac{3m_j}{4}.
\]
Moreover, since \(Y_k\le T_k^j\),
\[
\left\{
\frac1n\sum_{k=1}^nT_k^j<\frac{m_j}{2}
\right\}
\subset
\left\{
\frac1n\sum_{k=1}^nY_k<\frac{m_j}{2}
\right\}.
\]
Hence, by Hoeffding's inequality,
\[
\begin{split}
P\left(
\frac1n\sum_{k=1}^nT_k^j<\frac{m_j}{2}
\right)
&\le
P\left(
\frac1n\sum_{k=1}^nY_k<\frac{m_j}{2}
\right)\\
&\le
P\left(
\frac1n\sum_{k=1}^nY_k-E(Y_1)
<
-\frac{m_j}{4}
\right)\\
&\le
\exp\left\{
-\frac{2n(m_j/4)^2}{R^2}
\right\}.
\end{split}
\]
Therefore, there exist constants \(C_2,c_2>0\) such that
\[
P\left(
\frac1n\sum_{k=1}^nT_k^j<\frac{m_j}{2}
\right)
\le
C_2e^{-c_2n},
\qquad n\ge1.
\]
Therefore,
\[
\begin{split}
& P_i\left(
\frac{\sum_{k=1}^{t_j(n_t)}T_k^j}{t_j(n_t)}<\frac{m_j}{2},\,
\frac{t_j(n_t)}{n_t}\ge \pi_j-\epsilon,\,
n_t>d_jt
\right)
\\
&\qquad\le
\sum_{n\ge\lfloor(\pi_j-\epsilon)d_jt\rfloor}
P\left(
\frac1n\sum_{k=1}^nT_k^j<\frac{m_j}{2}
\right)
\le
C_3e^{-c_3t}
\end{split}
\]
for all sufficiently large \(t\).

On the complement of the two exceptional events above, if \(n_t>d_jt\), then
\[
\frac{t_j(n_t)}{n_t}\ge \pi_j-\epsilon>0
\]
and
\[
\frac{\sum_{k=1}^{t_j(n_t)}T_k^j}{t_j(n_t)}
\ge
\frac{m_j}{2}.
\]
Since \(r_j(t)\ge0\), we obtain
\[
\begin{split}
\frac{L_j(t)}{t}
&=
\frac{\sum_{k=1}^{t_j(n_t)}T_k^j+r_j(t)}{t_j(n_t)}
\frac{t_j(n_t)}{n_t}
\frac{n_t}{t}
\\
&\ge
\frac{m_j}{2}(\pi_j-\epsilon)\frac{n_t}{t}
>
\frac{m_j}{2}(\pi_j-\epsilon)d_j
=
c_{(j)}.
\end{split}
\]
Thus, on this event, \(L_j(t)>tc_{(j)}\). Consequently,
\[
P_i\big(L_j(t)\le tc_{(j)},\, n_t>d_jt\big)
\le
C_4e^{-c_4t}
\]
for suitable constants \(C_4,c_4>0\).

It remains to estimate \(P_i(n_t\le d_jt)\). Since \(d_j<1/\bar\mu\),
Lemma~\ref{lemma3.8} gives
\[
P_i(n_t\le d_jt)
\le
C_5t^{1-\eta}
\]
for all sufficiently large \(t\). Combining the last estimates, we obtain
\[
P_i\big(L_j(t)\le tc_{(j)}\big)
\le
C_4e^{-c_4t}+C_5t^{1-\eta}
\le
C_6t^{1-\eta}
\]
for all sufficiently large \(t\).

Since \(j\in\{1,\dots,K\}\) and the initial type \(i\in\{1,\dots,K\}\) were arbitrary,
and since \(K\) is finite, the constants can be chosen uniformly after increasing
\(C_6\). Therefore, for each \(j=1,\dots,K\), there exists \(c_{(j)}\in(0,1)\) such that
\[
P_i\left(
\frac{L_j(t)}{t}\le c_{(j)}
\right)
\le
C t^{-(\eta-1)},
\qquad t \text{ sufficiently large},\quad i=1,\dots,K.
\]
Hence Condition~\ref{Ac} holds with
\[
q=\eta-1.
\]
\(\hfill\square\)

\subsubsection*{Proof of Theorem \ref{Teorema2.6.8}:} By Proposition~\ref{Lema2.6.1}, for every \(\psi\in C_c^+(S)\),
\begin{equation}\label{eq:K-first-moment-bound}
\sup_{t\ge0}\sum_{i=1}^K
\int_{\mathbb R^N}
E_{(x,i)}\big(\langle X_t,\psi\rangle\big)\,dx
<\infty.
\end{equation}
We shall use this estimate with \(\psi=\phi\) and later with
\(\psi=C_\phi\bar f\).

Define, for \(\mathbf y=(y_1,\dots,y_K)\in[0,1]^K\),
\[
R_i(\mathbf y)
:=
\sum_{j=1}^K m_{i,j}y_j
-
\bigl(1-f_i(1-\mathbf y)\bigr),
\qquad i=1,\dots,K,
\]
where \(1-\mathbf y=(1-y_1,\dots,1-y_K)\). For an offspring vector
\(\mathbf n=(n_1,\dots,n_K)\in\mathbb N_0^K\), set
\[
A_{\mathbf n}(\mathbf y)
:=
\sum_{j=1}^K n_jy_j
-
1
+
\prod_{j=1}^K(1-y_j)^{n_j},
\qquad
\mathbf y\in[0,1]^K,
\]
with the convention \(a^0=1\) for \(a\in[0,1]\).

We first show that \(A_{\mathbf n}(\mathbf y)\ge0\). Let
\[
m:=n_1+\cdots+n_K.
\]
If \(m=0\), then \(A_{\mathbf n}(\mathbf y)=0\). Assume now that \(m\ge1\). Let
\(a_1,\dots,a_m\in[0,1]\) be the list obtained by repeating \(y_j\) exactly
\(n_j\) times, for \(j=1,\dots,K\). Then
\[
\prod_{j=1}^K(1-y_j)^{n_j}
=
\prod_{r=1}^m(1-a_r),
\qquad
\sum_{j=1}^K n_jy_j
=
\sum_{r=1}^m a_r.
\]
The elementary inequality
\[
1-\prod_{r=1}^m(1-a_r)\le \sum_{r=1}^m a_r,
\qquad 0\le a_r\le1,
\]
follows by induction on \(m\). Indeed, it is an equality for \(m=1\). If it
holds for \(m\), then, writing \(P_m:=\prod_{r=1}^m(1-a_r)\), we have
\[
1-P_m(1-a_{m+1})
=
1-P_m+P_m a_{m+1}
\le
\sum_{r=1}^m a_r+a_{m+1},
\]
because \(0\le P_m\le1\). Hence
\[
1-\prod_{j=1}^K(1-y_j)^{n_j}
\le
\sum_{j=1}^K n_jy_j.
\]
Therefore
\[
A_{\mathbf n}(\mathbf y)
=
\sum_{j=1}^K n_jy_j
-
\left[
1-\prod_{j=1}^K(1-y_j)^{n_j}
\right]
\ge0.
\]
Moreover, since
\[
0\le 1-\prod_{j=1}^K(1-y_j)^{n_j},
\]
we also have the upper bound
\[
0\le A_{\mathbf n}(\mathbf y)\le \sum_{j=1}^K n_jy_j.
\]

We next prove that \(A_{\mathbf n}\) is nondecreasing in each coordinate. Fix
\(\ell\in\{1,\dots,K\}\), and keep all coordinates except \(y_\ell\) fixed. If
\(n_\ell=0\), then \(A_{\mathbf n}\) does not depend on \(y_\ell\), and hence it is
constant in that coordinate. If \(n_\ell\ge1\), then, for \(0<y_\ell<1\),
\[
\frac{\partial A_{\mathbf n}}{\partial y_\ell}(\mathbf y)
=
n_\ell
-
n_\ell(1-y_\ell)^{n_\ell-1}
\prod_{j\ne\ell}(1-y_j)^{n_j}.
\]
Equivalently,
\[
\frac{\partial A_{\mathbf n}}{\partial y_\ell}(\mathbf y)
=
n_\ell
\left[
1-(1-y_\ell)^{n_\ell-1}
\prod_{j\ne\ell}(1-y_j)^{n_j}
\right].
\]
Since \(0\le 1-y_j\le1\) for every \(j\), we have
\[
0\le
(1-y_\ell)^{n_\ell-1}
\prod_{j\ne\ell}(1-y_j)^{n_j}
\le1.
\]
Thus
\[
\frac{\partial A_{\mathbf n}}{\partial y_\ell}(\mathbf y)\ge0,
\qquad 0<y_\ell<1.
\]
By the mean value theorem, \(A_{\mathbf n}\) is nondecreasing in \(y_\ell\) on
\((0,1)\). Since \(A_{\mathbf n}\) is continuous on \([0,1]^K\), the same
monotonicity extends to the boundary. Therefore \(A_{\mathbf n}\) is
nondecreasing in each coordinate on \([0,1]^K\).

We now average with respect to the offspring distribution of a type-\(i\)
particle. Let
\[
\boldsymbol\zeta_i=(\zeta_{i,1},\dots,\zeta_{i,K})
\]
be its offspring vector. By the definition of the offspring generating function,
\[
f_i(\mathbf s)
=
\mathbb E_i\left[
\prod_{j=1}^K s_j^{\zeta_{i,j}}
\right],
\qquad \mathbf s\in[0,1]^K.
\]
Taking \(\mathbf s=1-\mathbf y\), we get
\[
\mathbb E_i\left[
\prod_{j=1}^K(1-y_j)^{\zeta_{i,j}}
\right]
=
f_i(1-\mathbf y).
\]
Since
\[
\mathbb E_i(\zeta_{i,j})=m_{i,j},
\]
it follows that
\[
\begin{aligned}
\mathbb E_i\bigl(A_{\boldsymbol\zeta_i}(\mathbf y)\bigr)
&=
\mathbb E_i\left[
\sum_{j=1}^K \zeta_{i,j}y_j
-
1
+
\prod_{j=1}^K(1-y_j)^{\zeta_{i,j}}
\right]  \\
&=
\sum_{j=1}^K m_{i,j}y_j
-
1
+
f_i(1-\mathbf y) \\
&=
\sum_{j=1}^K m_{i,j}y_j
-
\bigl(1-f_i(1-\mathbf y)\bigr) \\
&=
R_i(\mathbf y).
\end{aligned}
\]
Furthermore, by the bound already proved,
\[
0\le
A_{\boldsymbol\zeta_i}(\mathbf y)
\le
\sum_{j=1}^K \zeta_{i,j}y_j
\le
\sum_{j=1}^K \zeta_{i,j}.
\]
Since \(\mathbb E_i(\zeta_{i,j})=m_{i,j}<\infty\) for each \(j\), the expectation
\(\mathbb E_i(A_{\boldsymbol\zeta_i}(\mathbf y))\) is finite. Consequently,
\[
R_i(\mathbf y)\ge0.
\]
Finally, if \(\mathbf y,\mathbf z\in[0,1]^K\) satisfy \(y_j\le z_j\) for every
\(j\), then the coordinatewise monotonicity of \(A_{\mathbf n}\) gives
\[
A_{\boldsymbol\zeta_i}(\mathbf y)
\le
A_{\boldsymbol\zeta_i}(\mathbf z)
\]
for every realization of \(\boldsymbol\zeta_i\). Taking expectations yields
\[
R_i(\mathbf y)\le R_i(\mathbf z).
\]
Thus \(R_i\) is nonnegative and nondecreasing in each coordinate on
\([0,1]^K\). Furthermore, for \(\mathbf y,\mathbf z\in[0,1]^K\), we claim that
\[
|A_{\mathbf n}(\mathbf y)-A_{\mathbf n}(\mathbf z)|
\le
2\sum_{j=1}^K n_j|y_j-z_j|.
\]
Indeed,
\[
\begin{aligned}
|A_{\mathbf n}(\mathbf y)-A_{\mathbf n}(\mathbf z)|
&\le
\sum_{j=1}^K n_j|y_j-z_j| \\
&\quad+
\left|
\prod_{j=1}^K(1-y_j)^{n_j}
-
\prod_{j=1}^K(1-z_j)^{n_j}
\right|.
\end{aligned}
\]
Writing the products as products of \(n_1+\cdots+n_K\) factors in \([0,1]\), the
elementary telescopic estimate
\[
\left|\prod_{r=1}^m a_r-\prod_{r=1}^m b_r\right|
\le
\sum_{r=1}^m |a_r-b_r|,
\qquad 0\le a_r,b_r\le1,
\]
gives
\[
\left|
\prod_{j=1}^K(1-y_j)^{n_j}
-
\prod_{j=1}^K(1-z_j)^{n_j}
\right|
\le
\sum_{j=1}^K n_j|y_j-z_j|.
\]
Hence
\[
|A_{\mathbf n}(\mathbf y)-A_{\mathbf n}(\mathbf z)|
\le
2\sum_{j=1}^K n_j|y_j-z_j|.
\]
Taking expectations with respect to the offspring distribution of a type-\(i\)
particle yields
\[
\begin{aligned}
|R_i(\mathbf y)-R_i(\mathbf z)|
&=
\left|
\mathbb E_i\left[
A_{\boldsymbol\zeta_i}(\mathbf y)
-
A_{\boldsymbol\zeta_i}(\mathbf z)
\right]
\right| \\
&\le
\mathbb E_i\left[
\left|
A_{\boldsymbol\zeta_i}(\mathbf y)
-
A_{\boldsymbol\zeta_i}(\mathbf z)
\right|
\right] \\
&\le
2\sum_{j=1}^K \mathbb E_i(\zeta_{i,j})|y_j-z_j| \\
&=
2\sum_{j=1}^K m_{i,j}|y_j-z_j|.
\end{aligned}
\]
Thus \(R_i\) is Lipschitz on \([0,1]^K\).

By the first-step decomposition for the Laplace functional,

\begin{equation}\label{V_t_equation}
\begin{split}
 V_t(x,i)
 &=
 (1-F_i(t))P_{t,\alpha_i}V_0(\cdot,i)(x)
 \\
 &\quad+
 \int_{[0,t]}
 P_{s,\alpha_i}
 \left(
 1-f_i\bigl(1-V_{t-s}(\cdot,1),\dots,1-V_{t-s}(\cdot,K)\bigr)
 \right)(x)\,F_i(ds).
\end{split}
\end{equation}

Integrating over \(\mathbb R^N\), using invariance of Lebesgue measure under
\(P_{t,\alpha_i}\), and using the definition of \(R_i\), we obtain
 \begin{equation}\label{eq:K-renewal-first-step}
 G_i(t)
 =
 (1-F_i(t))G_i(0)
 +
 \sum_{j=1}^K m_{i,j}\int_{[0,t]}G_j(t-s)\,F_i(ds)
 -
 \int_{[0,t]} B_i(t-s)\,F_i(ds),
 \end{equation}
where
\[
B_i(t):=
\int_{\mathbb R^N}
R_i\bigl(V_t(x,1),\dots,V_t(x,K)\bigr)\,dx.
\]

We prove that \(B_i\) is directly Riemann integrable. Let
\(D_\phi\subset\mathbb R^N\) be a ball such that
\[
\operatorname{supp}\phi\subset D_\phi\times\mathbf K.
\]
By Lemma~\ref{lem:finite-mean-local-bound}, applied to
\[
\mathcal A=D_\phi\times\mathbf K,
\]
there exists \(C_\phi>0\) such that,

\begin{equation*}
\begin{split}
V_t(x,i)
&=
E_{(x,i)}
\left(1-e^{-\langle X_t,\phi\rangle}\right)
\\
&\le
P_{(x,i)}
\left(X_t(\mathcal A)>0\right)
\le
C_{\phi}t^{-\rho},
\qquad t\ge1.
\end{split}
\end{equation*}

After increasing \(C_\phi\), this gives
\begin{equation}\label{eq:K-V-sup-bound}
\sup_{x\in\mathbb R^N,\ 1\le i\le K}V_t(x,i)
\le
C_\phi\min\{t^{-\rho},1\},
\qquad t\ge0.
\end{equation}
Also,
\[
0\le V_t(x,i)
\le
E_{(x,i)}\big(\langle X_t,\phi\rangle\big).
\]
Applying \eqref{eq:K-first-moment-bound} with \(\psi=\phi\), we get
\begin{equation}\label{eq:K-L1-bound}
\sup_{t\ge0}\sum_{i=1}^K\|V_t(\cdot,i)\|_1<\infty.
\end{equation}

We now use the branching assumption. Since \(u=\mathbf 1\) and \(M\) is stochastic,
 \[
 \sum_{j=1}^K m_{i,j}u_j=1,
 \qquad i=1,\ldots,K.
 \]
 Hence, for \(x\in[0,1]\),
 \[
 R_i(ux)
 =
 \sum_{j=1}^K m_{i,j}x-
 \bigl(1-f_i(1-ux)\bigr)
 =
 x-\bigl(1-f_i(1-ux)\bigr).
 \]
Therefore
\[
\begin{aligned}
\langle v,R(ux)\rangle
&=
\sum_{i=1}^K v_iR_i(ux)  \\
&=
\sum_{i=1}^K v_i
\left[
x-\bigl(1-f_i(1-ux)\bigr)
\right]  \\
&=
x\sum_{i=1}^K v_i
-
\sum_{i=1}^K v_i\bigl(1-f_i(1-ux)\bigr).
\end{aligned}
\]
Since \(\langle v,u\rangle=\sum_{i=1}^K v_i=1\), this gives
\[
\langle v,R(ux)\rangle
=
x-\langle v,1-\mathbf f(1-ux)\rangle.
\]
By the assumption on the branching mechanism,
\[
x-\langle v,1-\mathbf f(1-ux)\rangle
\sim
x^{1+\beta}L(x),
\qquad x\downarrow0.
\]
Consequently,
\[
\langle v,R(ux)\rangle
\sim
x^{1+\beta}L(x),
\qquad x\downarrow0.
\]

Since \(v_i>0\) for every \(i\) and \(R_i\ge0\), it follows that, for each fixed
\(i\),
\[
0\le v_iR_i(ux)
\le
\sum_{k=1}^K v_kR_k(ux)
=
\langle v,R(ux)\rangle.
\]
Therefore, for \(x>0\) sufficiently small,
\[
R_i(ux)
\le
\frac{1}{v_i}\langle v,R(ux)\rangle
\le
C_i x^{1+\beta}L(x).
\]
Since the type space is finite, increasing the constant if necessary, we may write
\[
R_i(ux)\le C x^{1+\beta}L(x),
\qquad i=1,\dots,K,
\]
for all sufficiently small \(x>0\).

Now let \(\mathbf y\in[0,1]^K\), and put
\[
x_{\mathbf y}:=\|\mathbf y\|_\infty.
\]
If \(x_{\mathbf y}=0\), then \(\mathbf y=0\) and \(R_i(\mathbf y)=R_i(0)=0\), so the
desired estimate is immediate. Assume \(x_{\mathbf y}>0\). Since \(u=\mathbf 1\),
we have the coordinatewise inequality
\[
\mathbf y\le ux_{\mathbf y}.
\]
Using the coordinatewise monotonicity of \(R_i\), we obtain
\[
R_i(\mathbf y)
\le
R_i(ux_{\mathbf y}).
\]
Thus, whenever \(x_{\mathbf y}=\|\mathbf y\|_\infty\) is sufficiently small,
\[
R_i(\mathbf y)
\le
C x_{\mathbf y}^{1+\beta}L(x_{\mathbf y})
=
C\|\mathbf y\|_\infty^{1+\beta}
L(\|\mathbf y\|_\infty).
\]

Finally, since \(L\) is slowly varying at zero, Potter's bound implies that, for
every \(\delta>0\), there exist constants \(C_\delta>0\) and \(r_\delta>0\) such
that
\[
L(r)\le C_\delta r^{-\delta},
\qquad 0<r<r_\delta.
\]
Taking \(r=\|\mathbf y\|_\infty\), we obtain, for
\(0<\|\mathbf y\|_\infty<r_\delta\),
\[
R_i(\mathbf y)
\le
C_\delta
\|\mathbf y\|_\infty^{1+\beta-\delta}.
\]
Hence, for every \(\delta\in(0,\beta)\), after possibly changing the constant,
\begin{equation}\label{eq:K-residual-bound}
R_i(\mathbf y)
\le
C_\delta
\|\mathbf y\|_\infty^{1+\beta-\delta},
\qquad
\|\mathbf y\|_\infty \text{ sufficiently small}.
\end{equation}

Since \(\rho>1/\beta\), we can choose \(\delta\in(0,\beta)\) sufficiently small so
that
\[
\rho(\beta-\delta)>1.
\]
Indeed, the inequality \(\rho>1/\beta\) is equivalent to \(\rho\beta>1\); hence,
by taking \(\delta>0\) small enough, we still have
\[
\rho\beta-\rho\delta>1.
\]
For all sufficiently large \(t\), \eqref{eq:K-V-sup-bound} implies that
\[
\max_{1\le j\le K}V_t(x,j)
\]
is uniformly small in \(x\). Hence, by \eqref{eq:K-residual-bound},
\[
\begin{split}
B_i(t)
&=
\int_{\mathbb R^N}
R_i\bigl(V_t(x,1),\dots,V_t(x,K)\bigr)\,dx
\\
&\le
C_\delta
\int_{\mathbb R^N}
\left(\max_{1\le j\le K}V_t(x,j)\right)^{1+\beta-\delta}\,dx
\\
&\le
C_\delta
\left(\sup_{x,j}V_t(x,j)\right)^{\beta-\delta}
\sum_{j=1}^K\|V_t(\cdot,j)\|_1
\\
&\le
C_\delta
\min\{t^{-\rho(\beta-\delta)},1\}.
\end{split}
\]
On bounded time intervals, the same domination follows from
\[
R_i(\mathbf y)\le \sum_{j=1}^K m_{i,j}y_j
\]
and \eqref{eq:K-L1-bound}. Thus, after increasing \(C_\delta\),
\begin{equation}\label{eq:K-Bi-dominated}
B_i(t)
\le
C_\delta
\min\{t^{-\rho(\beta-\delta)},1\},
\qquad t\ge0.
\end{equation}
The right-hand side is nonincreasing and integrable on \([0,\infty)\).

It remains to justify local Riemann integrability of \(B_i\).
Fix \(T>0\), and set
\[
\mathcal E:=L^1(\mathbb R^N)^K,
\qquad
\|\mathbf u\|_{\mathcal E}:=
\max_{1\le j\le K}\|u_j\|_1 .
\]
Since \(L^1(\mathbb R^N)\) is a Banach space, \(\mathcal E\), endowed with
\(\|\cdot\|_{\mathcal E}\), is also a Banach space.

We first prove that, for each \(i=1,\ldots,K\),
\[
t\mapsto V_t(\cdot,i)
\]
is càdlàg from \([0,T]\) into \(L^1(\mathbb R^N)\). Write the integral
equation \eqref{V_t_equation} as
\[
V_t(\cdot,i)
=
(1-F_i(t))P_{t,\alpha_i}V_0(\cdot,i)
+
\int_{[0,t]}
P_{s,\alpha_i}
\Phi_i\bigl(V_{t-s}(\cdot,1),\ldots,V_{t-s}(\cdot,K)\bigr)
\,F_i(ds),
\]
where
\[
\Phi_i(\mathbf y):=
1-f_i(1-y_1,\ldots,1-y_K),
\qquad
\mathbf y=(y_1,\ldots,y_K)\in[0,1]^K .
\]
For \(\mathbf y,\mathbf z\in[0,1]^K\),
\[
|\Phi_i(\mathbf y)-\Phi_i(\mathbf z)|
\le
\sum_{j=1}^K m_{i,j}|y_j-z_j|.
\]
Consequently, for functions \(\mathbf u,\mathbf v\) taking values in
\([0,1]^K\) a.e.,
\[
\|\Phi_i(\mathbf u)-\Phi_i(\mathbf v)\|_1
\le
\sum_{j=1}^K m_{i,j}\|u_j-v_j\|_1
\le
\|\mathbf u-\mathbf v\|_{\mathcal E},
\]
because \(M\) is stochastic.

The term
\[
t\mapsto (1-F_i(t))P_{t,\alpha_i}V_0(\cdot,i)
\]
is càdlàg in \(L^1(\mathbb R^N)\), since \(1-F_i\) is càdlàg and
\((P_{t,\alpha_i})_{t\ge0}\) is strongly continuous on \(L^1(\mathbb R^N)\).

Let \(\mathbf U\in D([0,T];\mathcal E)\), with \(0\le U_j(t,\cdot)\le1\)
a.e., and define
\[
I_i^{\mathbf U}(t)
:=
\int_{[0,t]}
P_{s,\alpha_i}
\Phi_i\bigl(U_1(t-s),\ldots,U_K(t-s)\bigr)
\,F_i(ds).
\]
We claim that \(t\mapsto I_i^{\mathbf U}(t)\) is càdlàg in
\(L^1(\mathbb R^N)\). Let \(t_n\downarrow t\). Then
\[
\begin{aligned}
I_i^{\mathbf U}(t_n)-I_i^{\mathbf U}(t)
&=
\int_{[0,t]}
P_{s,\alpha_i}
\left[
\Phi_i(\mathbf U(t_n-s))-\Phi_i(\mathbf U(t-s))
\right]
F_i(ds)
\\
&\quad+
\int_{(t,t_n]}
P_{s,\alpha_i}
\Phi_i(\mathbf U(t_n-s))
F_i(ds).
\end{aligned}
\]
The first term converges to zero in \(L^1\) by the right-continuity of
\(\mathbf U\), the Lipschitz bound for \(\Phi_i\), the contractivity of
\(P_{s,\alpha_i}\) on \(L^1\), and dominated convergence with respect to the
finite measure \(F_i\). For the second term, since
\(\mathbf U\in D([0,T];\mathcal E)\), it is bounded on \([0,T]\). Hence
\[
\sup_{0\le r\le T}
\|\Phi_i(\mathbf U(r))\|_1<\infty,
\]
and therefore
\[
\left\|
\int_{(t,t_n]}
P_{s,\alpha_i}
\Phi_i(\mathbf U(t_n-s))
F_i(ds)
\right\|_1
\le
C\,F_i((t,t_n]).
\]
Since \(F_i\) is right-continuous, \(F_i((t,t_n])\to0\). Thus
\(I_i^{\mathbf U}\) is right-continuous.

Similarly, if \(t_n\uparrow t\), then \(I_i^{\mathbf U}(t_n)\) converges in
\(L^1\) to
\[
\int_{[0,t)}
P_{s,\alpha_i}
\Phi_i\bigl(\mathbf U((t-s)-)\bigr)
F_i(ds),
\]
where \(\mathbf U(r-)\) denotes the left limit of \(\mathbf U\) at \(r>0\).
Indeed, for \(s<t\), \(t_n-s\uparrow t-s\), so
\[
\mathbf U(t_n-s)\longrightarrow \mathbf U((t-s)-)
\quad\text{in }\mathcal E,
\]
and dominated convergence applies with respect to the finite measure
\(F_i\) restricted to \([0,t)\). This proves the existence of left limits.
Hence \(I_i^{\mathbf U}\) is càdlàg.

Now choose \(h>0\) such that
\[
\max_{1\le i\le K}F_i(h)<1,
\]
which is possible because \(F_i(0)=0\) and \(F_i\) is right-continuous at
zero. Let \(\mathcal B_h\) denote the space of bounded Borel maps
\(\mathbf U:[0,h]\to\mathcal E\), endowed with the norm
\[
\|\mathbf U\|_{\infty,h}
:=
\sup_{0\le t\le h}\|\mathbf U(t)\|_{\mathcal E},
\]
and satisfying \(0\le U_j(t,\cdot)\le1\) a.e. for every \(j\) and \(t\).
This is a complete metric space. Indeed, since \(\mathcal E\) is a Banach
space, the space of all bounded maps from \([0,h]\) into \(\mathcal E\),
endowed with the supremum norm, is complete. The subspace of bounded Borel
maps is closed under uniform convergence, because uniform limits of Borel maps
with values in a metric space are Borel. Finally, the constraint
\(0\le U_j(t,\cdot)\le1\) a.e. is closed under \(L^1\)-convergence for each
fixed \(t\). Hence \(\mathcal B_h\) is complete. Moreover, the operator \(\mathcal T\) defined by the right-hand side of
\eqref{V_t_equation} maps \(\mathcal B_h\) into itself. Indeed, if
\(\mathbf U\in\mathcal B_h\), then \(\mathcal T\mathbf U\) is Borel
measurable. To see this, extend each \(U_j\) by zero outside \([0,h]\). Then
\[
(t,s)\mapsto U_j(t-s)\mathbf 1_{\{0\le s\le t\}}
\]
is Borel measurable from \([0,h]\times[0,h]\) into \(L^1(\mathbb R^N)\).
Since \(\Phi_i\) is Lipschitz as a map on \(L^1(\mathbb R^N)^K\), and
\((P_{s,\alpha_i})_{s\ge0}\) is strongly measurable on \(L^1(\mathbb R^N)\),
it follows that
\[
t\mapsto
\int_{[0,t]}
P_{s,\alpha_i}\Phi_i\bigl(\mathbf U(t-s)\bigr)\,F_i(ds)
\]
is Borel measurable as an \(L^1(\mathbb R^N)\)-valued map. Hence
\(\mathcal T\mathbf U\) is Borel measurable.

Furthermore, \(0\le V_0\le1\) and \(0\le \Phi_i(\mathbf y)\le1\) for
\(\mathbf y\in[0,1]^K\). Hence, for \(0\le t\le h\),
\[
0\le
(\mathcal T\mathbf U)_i(t)
\le
(1-F_i(t))P_{t,\alpha_i}\mathbf 1
+
\int_{[0,t]}P_{s,\alpha_i}\mathbf 1\,F_i(ds)
=
(1-F_i(t))+F_i(t)
=1
\qquad \text{a.e.}
\]
Moreover, for \(0\le t\le h\),
\[
\begin{aligned}
\|(\mathcal T\mathbf U)_i(t)\|_1
&\le
(1-F_i(t))\|P_{t,\alpha_i}V_0(\cdot,i)\|_1
+
\int_{[0,t]}
\|P_{s,\alpha_i}\Phi_i(\mathbf U(t-s))\|_1\,F_i(ds)
\\
&\le
\|V_0(\cdot,i)\|_1
+
\int_{[0,t]}
\|\Phi_i(\mathbf U(t-s))\|_1\,F_i(ds)
\\
&\le
\|V_0(\cdot,i)\|_1
+
\int_{[0,t]}
\sum_{j=1}^K m_{i,j}\|U_j(t-s)\|_1\,F_i(ds)
\\
&\le
\|V_0(\cdot,i)\|_1
+
F_i(h)\|\mathbf U\|_{\infty,h}
<\infty .
\end{aligned}
\]
Consequently, \(\mathcal T\mathbf U\in\mathcal B_h\).
Furthermore, for \(\mathbf U,\mathbf W\in\mathcal B_h\),
\[
\sup_{0\le t\le h}
\|(\mathcal T\mathbf U)(t)-(\mathcal T\mathbf W)(t)\|_{\mathcal E}
\le
\max_{1\le i\le K}F_i(h)
\sup_{0\le t\le h}
\|\mathbf U(t)-\mathbf W(t)\|_{\mathcal E}.
\]

Thus, by the contraction mapping theorem, the integral equation has a unique
bounded Borel solution on \([0,h]\). To see that this solution is càdlàg,
start the iteration from \(\mathbf U^{[0]}\equiv0\). Since \(\mathcal T\) maps
\(D([0,h];\mathcal E)\cap\mathcal B_h\) into itself, all iterates
\[
\mathbf U^{[n+1]}=\mathcal T\mathbf U^{[n]},\qquad n\ge0,
\]
are càdlàg. The contraction estimate implies that \(\mathbf U^{[n]}\) converges
uniformly in \(\mathcal E\) to the unique fixed point. Since the uniform limit
of càdlàg functions with values in a Banach space is càdlàg, the unique
solution is càdlàg. Moreover, by \(0\le V_t\le1\) and \eqref{eq:V-uniform-L1-bound}, the
probabilistic function
\[
\mathbf V(t):=(V_t(\cdot,1),\ldots,V_t(\cdot,K))
\]
belongs to \(\mathcal B_h\). It satisfies the same integral equation, and
therefore, by uniqueness in \(\mathcal B_h\), it coincides with the fixed
point. Hence
\[
t\mapsto V_t(\cdot,i)
\]
is càdlàg on \([0,h]\), for each \(i=1,\ldots,K\).

Repeating the argument on consecutive intervals gives the same conclusion on
\([0,T]\). Indeed, when the equation is solved on an interval \([a,a+h]\),
the unknown part of the convolution corresponds to \(0\le s\le t-a\le h\),
and is therefore controlled by the same contraction constant
\(\max_iF_i(h)<1\), while the remaining part only involves values already
constructed on earlier intervals. Therefore, for each \(i=1,\ldots,K\),
\[
t\mapsto V_t(\cdot,i)
\]
is càdlàg from \([0,T]\) into \(L^1(\mathbb R^N)\).

Since \(R_i\) is Lipschitz on \([0,1]^K\),
\[
|R_i(\mathbf y)-R_i(\mathbf z)|
\le
2\sum_{j=1}^K m_{i,j}|y_j-z_j|,
\qquad
\mathbf y,\mathbf z\in[0,1]^K,
\]
we have, for \(s,t\in[0,T]\),
\[
|B_i(t)-B_i(s)|
\le
C\sum_{j=1}^K
\|V_t(\cdot,j)-V_s(\cdot,j)\|_1 .
\]
Therefore \(B_i\) is càdlàg on compact intervals. Since every real-valued
càdlàg function on a compact interval is regulated, \(B_i\) is locally
Riemann integrable. In particular, \(B_i\) is Borel measurable and locally
bounded.

Together with \eqref{eq:K-Bi-dominated},
\[
0\le B_i(t)
\le
C_\delta \min\{t^{-\rho(\beta-\delta)},1\},
\qquad t\ge0,
\]
where \(\rho(\beta-\delta)>1\), this implies that \(B_i\) is directly Riemann
integrable on \([0,\infty)\). Indeed, \(B_i\) is locally Riemann integrable
and is dominated by an integrable nonincreasing function on \([0,\infty)\).

Define
\[
z_i(t)
:=
(1-F_i(t))G_i(0)
-
\int_{[0,t]}B_i(t-s)\,F_i(ds).
\]
The function \(1-F_i\) is nonincreasing and integrable on \([0,\infty)\),
since
\[
\int_0^\infty (1-F_i(t))\,dt=m_i<\infty.
\]
Hence \(1-F_i\) is directly Riemann integrable. Since \(B_i\) is nonnegative
and directly Riemann integrable, and \(F_i\) is a finite measure on
\([0,\infty)\), the convolution
\[
t\mapsto
\int_{[0,t]}B_i(t-s)\,F_i(ds)
\]
is directly Riemann integrable; see \cite[Theorem~1]{Sgibnev2023}.
Therefore \(z_i\), being the difference of two directly Riemann integrable
functions, is directly Riemann integrable in the real-valued sense.

Moreover, \(z_i\) is locally bounded. Indeed, for every \(T>0\),
\[
0\le (1-F_i(t))G_i(0)\le G_i(0),
\qquad 0\le t\le T.
\]
Since \(B_i\) is càdlàg, it is bounded on compact intervals. Therefore
\[
\sup_{0\le t\le T}
\int_{[0,t]}B_i(t-s)\,F_i(ds)
\le
\sup_{0\le r\le T}B_i(r)\,F_i([0,T])
\le
\sup_{0\le r\le T}B_i(r)
<\infty .
\]
Thus \(z_i\) is Borel measurable, locally bounded, and directly Riemann
integrable.

With this notation, \eqref{eq:K-renewal-first-step} becomes
\begin{equation}\label{eq:K-renewal-final-system}
 G_i(t)
 =
 z_i(t)
 +
 \sum_{j=1}^K
 \int_{[0,t]}G_j(t-s)\,F_{i,j}(ds),
 \end{equation}
where
\[
F_{i,j}(t):=m_{i,j}F_i(t).
\]
As in Proposition~\ref{Lema2.6.1}, the matrix of measures
\((F_{i,j})_{i,j=1}^K\) satisfies
\[
F(\infty)=M,
\qquad
p(F(\infty))=1,
\qquad
F(0)=0,
\qquad
p(F(0))=0<1,
\]
and is non-lattice because \(m_{i,i}>0\) and \(F_i\) is non-arithmetic. Since each
\(z_i\) is directly Riemann integrable, Theorem~2.2\textup{(iii)} of \cite{Athreya2}
applies to \eqref{eq:K-renewal-final-system}. Thus, for every \(i=1,\dots,K\),
\begin{equation}\label{eq:K-Gi-limit}
\lim_{t\to\infty}G_i(t)
=
\frac{
\sum_{k=1}^K\pi_k\int_0^\infty z_k(s)\,ds
}{\bar m}.
\end{equation}

By Tonelli's theorem,
\[
\begin{split}
 \int_0^\infty
 \int_{[0,t]}B_i(t-s)\,F_i(ds)\,dt
 &=
 \int_{[0,\infty)}
 \int_s^\infty B_i(t-s)\,dt\,F_i(ds)
 \\
 &=
 \int_{[0,\infty)}
 \int_0^\infty B_i(r)\,dr\,F_i(ds)
 =
 \int_0^\infty B_i(r)\,dr.
 \end{split}
 \]
Therefore,
\[
\int_0^\infty z_i(t)\,dt
=
m_iG_i(0)
-
\int_0^\infty B_i(s)\,ds.
\]
Substituting this into \eqref{eq:K-Gi-limit}, one obtains
\[
\lim_{t\to\infty}G_i(t)
=
\frac{
\sum_{k=1}^K\pi_km_kG_k(0)
-
\sum_{k=1}^K\pi_k\int_0^\infty B_k(s)\,ds
}{\bar m}.
\]
Since
\[
\theta_k=\frac{\pi_km_k}{\bar m},
\qquad
\sum_{k=1}^K\theta_k=1,
\]
we conclude that
\begin{equation}\label{eq:K-G-limit}
G(\phi)
:=
\lim_{t\to\infty}g_\phi(t)
=
\sum_{k=1}^K\theta_kG_k(0)
-
\frac1{\bar m}
\sum_{k=1}^K\pi_k\int_0^\infty B_k(s)\,ds.
\end{equation}

It remains to prove the first-order identity and the positivity of \(G(\phi)\). Let
\[
\bar f_i(x):=1-e^{-\phi_i(x)},
\qquad i=1,\dots,K.
\]
For \(0<\epsilon<1/2\), define
\[
\phi_{\epsilon,i}(x):=-\log(1-\epsilon\bar f_i(x)),
\qquad
\phi_\epsilon(x,i):=\phi_{\epsilon,i}(x),
\]
and set
\[
V_t^\epsilon(x,i)
:=
E_{(x,i)}
\left(
1-e^{-\langle X_t,\phi_\epsilon\rangle}
\right).
\]
Then
\[
V_0^\epsilon(x,i)=\epsilon\bar f_i(x).
\]
Moreover, \(0\le \phi_\epsilon\le\phi\), and hence
\[
0\le V_t^\epsilon(x,i)\le V_t(x,i).
\]
Since \(x\mapsto-\log(1-x)/x\) is bounded on \([0,1/2]\), there exists
\(C_\phi>0\) such that
\[
\frac{\phi_{\epsilon,i}(x)}{\epsilon}
\le
C_\phi\bar f_i(x).
\]
Using \(1-e^{-y}\le y\), we obtain
\begin{equation}\label{eq:K-Veps-linear-bound}
\frac{V_t^\epsilon(x,i)}{\epsilon}
\le
E_{(x,i)}
\left(
\left\langle X_t,C_\phi\bar f\right\rangle
\right).
\end{equation}

Let
\[
B_i^\epsilon(t)
:=
\int_{\mathbb R^N}
R_i\bigl(V_t^\epsilon(x,1),\dots,V_t^\epsilon(x,K)\bigr)\,dx.
\]
We prove that
\begin{equation}\label{eq:K-Beps-zero}
\lim_{\epsilon\downarrow0}
\frac1\epsilon
\int_0^\infty B_i^\epsilon(s)\,ds
=
0.
\end{equation}
For all sufficiently large \(t\), using \eqref{eq:K-residual-bound},
\(V_t^\epsilon\le V_t\), and \eqref{eq:K-Veps-linear-bound},
\[
\begin{split}
\frac{B_i^\epsilon(t)}{\epsilon}
&\le
C
\int_{\mathbb R^N}
\left(\max_{1\le j\le K}V_t(x,j)\right)^{\beta-\delta}
\sum_{j=1}^K\frac{V_t^\epsilon(x,j)}{\epsilon}\,dx
\\
&\le
C
\left(\sup_{x,j}V_t(x,j)\right)^{\beta-\delta}
\sum_{j=1}^K
\int_{\mathbb R^N}
E_{(x,j)}
\left(
\left\langle X_t,C_\phi\bar f\right\rangle
\right)\,dx
\\
&\le
C\min\{t^{-\rho(\beta-\delta)},1\}.
\end{split}
\]
The right-hand side is integrable. On bounded time intervals, domination follows from
\[
R_i(\mathbf y)\le\sum_{j=1}^K m_{i,j}y_j,
\]
from \eqref{eq:K-Veps-linear-bound}, and from \eqref{eq:K-first-moment-bound}. Thus
\(B_i^\epsilon(t)/\epsilon\) is dominated by an integrable function independent of
\(\epsilon\).

For each fixed \(t\ge0\), we now show that
\[
\frac{B_i^\epsilon(t)}{\epsilon}\longrightarrow0.
\]
For each fixed \(x\in\mathbb R^N\) and \(j\in\{1,\dots,K\}\),
\[
V_t^\epsilon(x,j)
=
E_{(x,j)}
\left(
1-e^{-\langle X_t,\phi_\epsilon\rangle}
\right)
\longrightarrow0,
\qquad \epsilon\downarrow0,
\]
by dominated convergence. In addition, \eqref{eq:K-Veps-linear-bound} gives
\[
V_t^\epsilon(x,j)
\le
\epsilon
E_{(x,j)}
\left(
\left\langle X_t,C_\phi\bar f\right\rangle
\right).
\]
Hence, for fixed \(t,x,j\),
\[
V_t^\epsilon(x,j)=O(\epsilon),
\qquad \epsilon\downarrow0.
\]
Using \eqref{eq:K-residual-bound}, we obtain
\[
\frac{
R_i\bigl(V_t^\epsilon(x,1),\dots,V_t^\epsilon(x,K)\bigr)
}{\epsilon}
\longrightarrow0.
\]
Furthermore,
\[
\frac{R_i(V_t^\epsilon(x,1),\dots,V_t^\epsilon(x,K))}{\epsilon}
\le
\sum_{j=1}^K m_{i,j}
\frac{V_t^\epsilon(x,j)}{\epsilon}
\le
\sum_{j=1}^K m_{i,j}
E_{(x,j)}
\left(
\left\langle X_t,C_\phi\bar f\right\rangle
\right),
\]
and the right-hand side is integrable in \(x\) by \eqref{eq:K-first-moment-bound}, applied
with \(\psi=C_\phi\bar f\). Therefore, dominated convergence in \(x\) gives
\[
\frac{B_i^\epsilon(t)}{\epsilon}\longrightarrow0
\qquad
\text{for every fixed }t\ge0.
\]
A second application of dominated convergence in \(t\), using the integrable domination
obtained above, proves \eqref{eq:K-Beps-zero}.

Applying \eqref{eq:K-G-limit} with \(\phi\) replaced by \(\phi_\epsilon\), we obtain
\[
G(\phi_\epsilon)
=
\sum_{k=1}^K\theta_k\|V_0^\epsilon(\cdot,k)\|_1
-
\frac1{\bar m}
\sum_{k=1}^K\pi_k\int_0^\infty B_k^\epsilon(s)\,ds.
\]
Since
\[
\|V_0^\epsilon(\cdot,k)\|_1
=
\epsilon\|\bar f_k\|_1,
\]
division by \(\epsilon\) gives
\[
\frac{G(\phi_\epsilon)}{\epsilon}
=
\sum_{k=1}^K\theta_k\|\bar f_k\|_1
-
\frac1{\bar m}
\sum_{k=1}^K\pi_k
\frac1\epsilon
\int_0^\infty B_k^\epsilon(s)\,ds.
\]
Using \eqref{eq:K-Beps-zero}, we obtain
\[
\lim_{\epsilon\downarrow0}
\frac{G(\phi_\epsilon)}{\epsilon}
=
\sum_{k=1}^K\theta_k\|\bar f_k\|_1.
\]
Since
\[
V_0(x,k)=1-e^{-\phi_k(x)}=\bar f_k(x),
\]
the right-hand side equals
\[
\sum_{k=1}^K\theta_k\|V_0(\cdot,k)\|_1.
\]
This proves \eqref{eq:K-first-order-limit}.

Finally, if \(\phi\neq0\), then
\[
\sum_{k=1}^K\theta_k\|\bar f_k\|_1>0,
\]
because \(\theta_k>0\) for every \(k\). Hence
\[
G(\phi_\epsilon)>0
\]
for all sufficiently small \(\epsilon>0\). Since \(0\le\phi_\epsilon\le\phi\), monotonicity gives
\[
0\le V_t^\epsilon(x,i)\le V_t(x,i),
\qquad
t\ge0,\ x\in\mathbb R^N,\ i=1,\dots,K.
\]
Therefore
\[
g_{\phi_\epsilon}(t)\le g_\phi(t),
\qquad t\ge0.
\]
Letting \(t\to\infty\), and using the existence of both limits, yields
\[
G(\phi_\epsilon)\le G(\phi).
\]
Consequently,
\[
G(\phi)>0.
\]

\(\hfill\square\)

\subsection*{Auxiliary Results}

\begin{lemma}\label{lem:finite-mean-local-bound}
Assume that Condition~\ref{Ac} holds and set
\[
\alpha_*:=\alpha_1\wedge\cdots\wedge\alpha_K,
\qquad
\rho:=q\wedge\frac{N}{\alpha_*}.
\]
Then, for every \(\phi\in C_c^+(S)\), there exists \(C_\phi>0\) such that, for all \(t\ge
1\) and all \((x,i)\in S\),
\[
V_t(x,i)
\le
P_{(x,i)}\bigl(X_t(D_\phi\times\mathbf K)>0\bigr)
\le
C_\phi t^{-\rho},
\]
where \(D_\phi\subset\mathbb R^N\) is any ball such that
\(\operatorname{supp}\phi\subset D_\phi\times\mathbf K\).
\end{lemma}

\begin{proof}
Following the argument of Proposition~4.2 of \cite{GRW}, and using
Condition~\ref{Ac}, there exists \(C_\phi>0\) such that, for all \(t\ge1\) and all
\((x,i)\in S\),
\[
P_{(x,i)}\bigl(X_t(D_\phi\times\mathbf K)>0\bigr)
\le
C_\phi\left(t^{-q}+t^{-N/\alpha_*}\right)
\le
C_\phi t^{-\rho},
\]
where \(\alpha_*=\alpha_1\wedge\cdots\wedge\alpha_K\) and
\(\rho=q\wedge N/\alpha_*\). Since
\[
1-e^{-\langle X_t,\phi\rangle}
\le
1_{\{X_t(D_\phi\times\mathbf K)>0\}},
\]
we obtain
\[
V_t(x,i)
=
E_{(x,i)}\left(1-e^{-\langle X_t,\phi\rangle}\right)
\le
P_{(x,i)}\bigl(X_t(D_\phi\times\mathbf K)>0\bigr)
\le
C_\phi t^{-\rho}.
\]
This proves the lemma.
\end{proof}

The following lemmas provide auxiliary results that will be used in the proof of Theorem~\ref{Teo2.3}.

\begin{lemma}\label{lemaaux1}
Let \(T^{(1)},\dots,T^{(K)}\) be independent nonnegative random variables, where
\(T^{(i)}\) has distribution \(F_i\). Assume that
\[
1-F_i(t)\le A t^{-\eta^{(i)}},
\qquad t>0,\quad i=1,\dots,K,
\]
for some \(A>0\) and \(\eta^{(i)}>1\). Set
\[
\eta:=\min_{1\le i\le K}\eta^{(i)}
\]
and
\[
T^*:=\max_{1\le i\le K}T^{(i)}.
\]
Then there exists a constant \(\bar A>0\) such that
\begin{equation}\label{eq:max-tail-K}
P(T^*>t)\le \bar A t^{-\eta},
\qquad t>0.
\end{equation}
In particular, \(E(T^*)<\infty\).
\end{lemma}

\begin{proof}
By the union bound,
\[
P(T^*>t)
=
P\left(\bigcup_{i=1}^K\{T^{(i)}>t\}\right)
\le
\sum_{i=1}^K P(T^{(i)}>t).
\]
For \(t\ge1\), since \(\eta\le\eta^{(i)}\) for every \(i\),
\[
P(T^*>t)
\le
\sum_{i=1}^K A t^{-\eta^{(i)}}
\le
KA t^{-\eta}.
\]
For \(0<t<1\), we have \(P(T^*>t)\le1\le t^{-\eta}\). Hence
\[
P(T^*>t)\le \bar A t^{-\eta},
\qquad t>0,
\]
with \(\bar A:=KA\vee1\). Since \(\eta>1\),
\[
E(T^*)
=
\int_0^\infty P(T^*>t)\,dt
<\infty.
\]
This proves the lemma.
\end{proof}

\begin{lemma}\label{lemaaux2}
Let \(Y_1,Y_2,\dots\) be i.i.d. nonnegative random variables with mean
\(m<\infty\). Assume that
\[
P(Y_1>t)\le A t^{-\eta},
\qquad t>0,
\]
for some \(A>0\) and \(\eta>1\). Then, for every \(c>0\), there exists
\(C>0\) such that, for all sufficiently large \(n\) and all \(x\ge cn\),
\begin{equation}\label{cotaeta1}
P\left(\sum_{k=1}^nY_k-nm\ge x\right)
\le C n x^{-\eta}.
\end{equation}
Consequently, for every \(\delta>0\), there exists \(C_\delta>0\) such that,
for all sufficiently large \(n\),
\begin{equation}\label{cotaeta2}
P\left(\frac1n\sum_{k=1}^nY_k-m\ge\delta\right)
\le C_\delta n^{1-\eta}.
\end{equation}
\end{lemma}

\begin{proof}
Apply Theorem~2 of \cite{Naga} to the centered variables
\[
X_k:=Y_k-m,
\qquad k\ge1.
\]
Then \(E(X_1)=0\). Moreover, for every \(y>0\),
\[
P(X_1>y)
=
P(Y_1>y+m)
\le
A(y+m)^{-\eta}
\le
A y^{-\eta}.
\]
Therefore, Theorem~2 of \cite{Naga} gives, for each fixed \(c>0\), a constant
\(C>0\) such that, for all sufficiently large \(n\) and all \(x\ge cn\),
\[
P\left(\sum_{k=1}^nX_k\ge x\right)
\le
C n x^{-\eta}.
\]
Since
\[
\sum_{k=1}^nX_k
=
\sum_{k=1}^nY_k-nm,
\]
this proves \eqref{cotaeta1}.

Taking \(x=\delta n\) in \eqref{cotaeta1}, we obtain
\[
P\left(
\frac1n\sum_{k=1}^nY_k-m\ge\delta
\right)
=
P\left(
\sum_{k=1}^nY_k-nm\ge \delta n
\right)
\le
C n(\delta n)^{-\eta}.
\]
Hence
\[
P\left(
\frac1n\sum_{k=1}^nY_k-m\ge\delta
\right)
\le
C_\delta n^{1-\eta},
\]
where
\[
C_\delta:=C\delta^{-\eta}.
\]
This proves \eqref{cotaeta2}.
\end{proof}

\begin{lemma}\label{lemma3.8}
Assume \eqref{tail2}. Let \(Z=(Z_t)_{t\ge0}\) be the Markov renewal process on
\(\mathbf K=\{1,\dots,K\}\), and let \(n_t\) be the number of renewals up to time \(t\).
Let \(T^{(1)},\dots,T^{(K)}\) be independent random variables such that \(T^{(r)}\) has
distribution \(F_r\), and define
\[
T^*:=\max_{1\le r\le K}T^{(r)}.
\]
Set
\[
\bar\mu:=E(T^*)<\infty.
\]
If \(c>0\) satisfies
\[
c<\frac1{\bar\mu},
\]
then, for every initial type \(i\in\mathbf K\), there exists a constant \(C>0\) such that,
for all sufficiently large \(t\),
\begin{equation}\label{cotan}
P_i(n_t\le ct)
\le
C t^{1-\eta},
\end{equation}
where
\[
\eta:=\min_{1\le r\le K}\eta^{(r)}.
\]
\end{lemma}

\begin{proof}
Let \((Y_n)_{n\ge0}\) be the embedded Markov chain of the Markov renewal
process. Thus \(Y_0\) is the initial type. Conditionally on
\((Y_n)_{n\ge0}\), let the holding times be independent random variables
\((\widetilde T_n)_{n\ge0}\), where \(\widetilde T_n\) has distribution
\(F_{Y_n}\). Define
\[
\sigma_0:=0,
\qquad
\sigma_n:=\sum_{r=0}^{n-1}\widetilde T_r,
\qquad n\ge1.
\]
Then
\[
n_t=\max\{n\ge0:\sigma_n\le t\}.
\]
Put
\[
n:=\lfloor ct\rfloor+1.
\]
Then
\[
\{n_t\le ct\}\subseteq \{\sigma_n>t\}.
\]

We now compare the holding times of the Markov renewal process with
independent copies of \(T^*\). Let \(F_*\) be the distribution function of
\[
T^*:=\max_{1\le r\le K}T^{(r)}.
\]
Since \(T^{(r)}\le_{\mathrm{st}}T^*\) for every \(r=1,\ldots,K\), we have
\[
F_r(x)\ge F_*(x),
\qquad x\ge0,\quad r=1,\ldots,K.
\]
Let \(U_0,U_1,\ldots\) be i.i.d. random variables uniformly distributed on
\((0,1)\), independent of the embedded chain \((Y_n)_{n\ge0}\). Using the
usual generalized inverse
\[
F^{-1}(u):=\inf\{x\ge0:F(x)\ge u\},
\qquad 0<u<1,
\]
define
\[
\widetilde T_n:=F_{Y_n}^{-1}(U_n),
\qquad
T_n^*:=F_*^{-1}(U_n),
\qquad n\ge0.
\]
Then, conditionally on \((Y_n)_{n\ge0}\), the variables
\((\widetilde T_n)_{n\ge0}\) are independent and \(\widetilde T_n\) has
distribution \(F_{Y_n}\). Hence this construction has the law of the holding
times of the Markov renewal process. Moreover, since
\(F_{Y_n}\ge F_*\), the monotonicity of generalized inverses gives
\[
\widetilde T_n\le T_n^*,
\qquad n\ge0,
\]
almost surely. The variables \(T_0^*,T_1^*,\ldots\) are i.i.d. with
distribution \(F_*\), that is, with the distribution of \(T^*\). Consequently,
for every \(n\ge1\),
\[
\sigma_n
=
\sum_{r=0}^{n-1}\widetilde T_r
\le
\sum_{r=0}^{n-1}T_r^*
=:\bar S_n
\qquad\text{a.s.}
\]
Therefore,
\[
P_i(n_t\le ct)
\le
P_i(\sigma_n>t)
\le
P(\bar S_n>t).
\]

Since \(n=\lfloor ct\rfloor+1\), we have \(n/t\to c\). Because
\(c<1/\bar\mu\),
\[
\frac{t}{n}\longrightarrow \frac1c>\bar\mu.
\]
Therefore, there exists \(\delta>0\) such that, for all sufficiently large
\(t\),
\[
t-n\bar\mu\ge \delta n.
\]
Thus,
\[
P(\bar S_n>t)
=
P(\bar S_n-n\bar\mu>t-n\bar\mu)
\le
P(\bar S_n-n\bar\mu\ge \delta n).
\]

By Lemma~\ref{lemaaux1}, the distribution of \(T^*\) satisfies
\[
P(T^*>x)\le \bar A x^{-\eta},
\qquad x>0.
\]
 Applying Lemma~\ref{lemaaux2} to the i.i.d. sequence
 \[
 Y_k:=T_{k-1}^*,
 \qquad k\ge1,
 \]
 with mean \(\bar\mu\), we obtain
 \[
 P(\bar S_n-n\bar\mu\ge \delta n)
 \le
 C n^{1-\eta}
 \]
 for all sufficiently large \(n\). Since \(n=\lfloor ct\rfloor+1\), it follows
that
\[
P(\bar S_n-n\bar\mu\ge \delta n)
\le
C' t^{1-\eta}
\]
for all sufficiently large \(t\). Combining the preceding estimates gives
\[
P_i(n_t\le ct)
\le
C' t^{1-\eta},
\]
which proves \eqref{cotan}.
\end{proof}
\section*{Acknowledgements}
J.H.R.G. is  supported by FAPESP grants  2025/03804-0.

\bibliographystyle{alea3}

\bibliography{ALEA504}

\end{document}